\documentclass[12pt]{article}
\usepackage{a4wide}

\usepackage{latexsym,amssymb,amsmath,epsfig,amsfonts,algorithm}
\usepackage{multirow,multicol,array,bm}
\usepackage[authoryear]{natbib}
\usepackage[colorlinks=true,allcolors=blue]{hyperref}
\usepackage{float}
\usepackage{caption}
\usepackage[utf8]{inputenc}
\usepackage[english]{babel}
\usepackage[noend]{algpseudocode}
\usepackage{tikz}
\usetikzlibrary{arrows}
\usepackage{bbm}
\usepackage[noend]{algpseudocode}
\usepackage{graphicx}
\usepackage[skip=0.333\baselineskip]{subcaption}
\usepackage[font=large,labelfont=bf]{caption}
\usepackage{mathtools}

\newcommand{\ds}{\displaystyle}

\newcommand*{\Scale}[2][4]{\scalebox{#1}{$#2$}}%
\providecommand{\keywords}[1]{\textbf{\textit{Keywords---}} #1}
\newtheorem{definition}{Definition}[section]  
\newtheorem{theorem}{Theorem}[section] 
\newtheorem{lemma}{Lemma}[section]

\newtheorem{proposition}{Proposition}[section] 
\newtheorem{paradox}{Paradox}[section] 
\newtheorem{corollary}{Corollary}[section] 
\numberwithin{figure}{section}
\numberwithin{table}{section}
\numberwithin{equation}{section}
\begin{document}

\title{Strategic Customer Behavior in an $M/M/1$ Feedback Queue}


\author{Mark Fackrell \and  Peter Taylor \and Jiesen Wang        
}

\author{Mark Fackrell, Peter Taylor, Jiesen Wang\\
	\normalsize{School of Mathematics and Statistics, The University of Melbourne,}\\
	\normalsize{Victoria 3010, Australia}\\
	\normalsize{contact:\{fackrell, taylorpg\}@unimelb.edu.au, jiesenw@student.unimelb.edu.au}
}
\date{February 12, 2021}

\maketitle

\begin{abstract}
We investigate the behavior of equilibria in an $M/M/1$ feedback queue where price and time sensitive customers are homogeneous with respect to service valuation and cost per unit time of waiting. Upon arrival, customers can observe the number of customers in the system and then decide to join or to balk. Customers are served in order of arrival. After being served, each customer either successfully completes the service and departs the system with probability $q$, or the service fails and the customer immediately joins the end of the queue to wait to be served again until she successfully completes it. 

We analyse this decision problem as a noncooperative game among the customers. We show that there exists a unique symmetric Nash equilibrium threshold strategy. 
We then prove that the symmetric Nash equilibrium threshold strategy is evolutionarily stable. Moreover, if we relax the strategy restrictions by allowing customers to renege, in the new Nash equilibrium, customers have a greater incentive to join. However, this does not necessarily increase the equilibrium expected payoff, and for some parameter values, it decreases it.

\end{abstract}

\section{Introduction} \label{sec:intro}

Consider an $M/M/1$ first-come-first-served queue where customers arrive according to a Poisson process with rate $\lambda$ and the service times for each customer are independently and identically distributed according to an exponential distribution with parameter $\mu$. After being served, each customer either successfully completes the service and departs from the system with probability $q$, or the service fails and the customer immediately joins the end of the queue to wait to be served again until she successfully completes it.

We define the sojourn time as the total time a customer spends in the system, so it includes both the waiting time and the service time. Upon arriving at the queue, the newly arrived customer observes the number of customers in the system, and by considering the trade-off between her expected sojourn time and the reward due to a successful service completion, she makes a decision to join the queue or balk depending on the number of customers present when she arrives. 

\begin{figure}
	\centering
	\includegraphics[width=10cm]{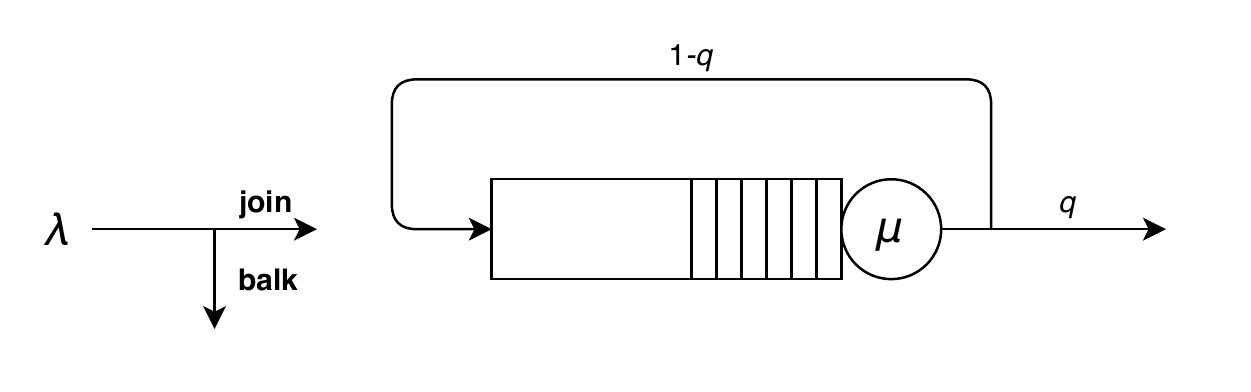}
	\caption{An $M/M/1$ feedback queue with strategic customers when reneging is not allowed.} \label{F1}
\end{figure}

The cost is assumed to be linear in the sojourn time with rate $C$. To non-dimensionalise the model, we set $C=1$ for the rest of the paper. The reward to a customer when she successfully completes her service, which is assumed to be identical across customers, is denoted by $R_0$. Let $R$ be the reward that a customer actually obtains when she leaves the system. In this first model, customers are not allowed to leave until they successfully complete their service. Hence, the random variable $R$ is equal to the reward $R_0$ with probability one, but the reason that we have introduced it is that the reward is truly random for the system with reneging that we consider later. Indeed, for that system the random variable $R$ is equal to the reward $R_0$ with some probability less than one and equal to zero with some positive probability. Customers decide to join as long as their expected payoff, which is defined as the difference between their expected reward and their expected cost, is positive. See Figure \ref{F1} for an illustration of the system. The sojourn time of a customer depends on the service times of all customers that are served before she leaves the system and, if she has to repeat her service, it is possible that some of these services are for customers who joined the queue after her. It follows that her expected reward depends on the joining strategy of other customers. As a consequence, the best response of each customer is a function of both the position at which she joins the system and the other customers’ strategies. For this reason, it is natural to consider the decision problem in a game theoretic framework and to look into the Nash equilibrium strategy for each customer (see \citet{HH03}). 

The study of the instantaneous Bernoulli feedback queue goes back to \citet{T63}, in which he obtained the expected total waiting time for the $M/G/1$ feedback queue by deriving the joint transform of the distribution of the queue length and the remaining service time. \citet {DMS80}, \citet{D81}, and \citet{DK84} further studied the queue length, the total sojourn time, and the waiting time. 
\citet{T87} applied the instantaneous Bernoulli feedback queue to study packet transmissions in an error-prone channel with probability of successful transmission $q$. 
This transmission style is similar to segmented message transmission with the number of segments in a message geometrically distributed with mean $1/q$. \citet{FH05}, and \citet{GH09} studied queueing networks with feedback loops and intelligent customers. However, the customers in their model are not strategic in the sense that their decisions do not depend on the others' decisions, which is different from our setting. 

\citet{AS98} analysed a system of observable egalitarian processor sharing queues, where customers decide to join or to balk after observing the number of customers in the system upon their arrival, and are not allowed to renege at any stage after joining. They calculated numerically the symmetric threshold equilibrium strategy for the case of Poisson arrivals and proposed a dynamic learning scheme, which converges to the symmetric Nash equilibrium strategy. 
\citet{AKLL18} considered a polling system with two queues where a single server serves the two nodes in a cyclic fashion with exhaustive service. Customers can choose which queue to join upon arrival. They analysed the Nash equilibrium strategies under three scenarios of available information of the queue lengths or the position of the server at decision epochs, and obtained the Nash equilibrium strategies via a new iterative algorithm.
In both \cite{AS98} and \cite{AKLL18}, customers' best response is affected by future arrivals.

In this paper, we study an $M/M/1$ queue with instantaneous Bernoulli feedback, and allow each customer to determine whether to join the queue or not after observing the number of customers in the system. Similar to \citet{AS98}, in our model, a tagged customer's sojourn time is affected by the joining behavior of future arrivals.
This model was first analysed in an unpublished technical report by \citet{BC13}. They considered a first-come-first-served $GI/G/1$ Bernoulli feedback queue with arriving customers observing the number of customers in the system before deciding whether to join or not, but not allowed to renege. 

Although arriving customers see the stationary distribution of the number of customers in the system, due to the balking and different joining positions, the distribution observed by feedback customers requires further analysis (see \citet[Section 2.10]{W88}, \citet{BVD97}), which makes the expected sojourn time computation nontrivial. In this paper, we efficiently obtain the expected payoff of a joining customer for any parameter set using matrix analytic methods (see \citet{N81}), which can also be easily extended to other models. In particular, we compute a customer's conditional expected payoff based on her joining position and the other customers’ threshold values, by solving Poisson's equation for a discrete-time nonhomogeneous quasi-birth-and-death-process. Then we explicitly propose the Nash equilibrium strategies (pure or mixed) of threshold type. 

Every time a customer joins at the end of the queue due to a service failure, the time she has already spent becomes a sunk cost. Also, it is possible that her conditions have deteriorated with time. For example, the system could have been empty when a tagged customer first arrived at the queue, but has become overcrowded before she goes to the end of the queue due to a service failure, because of a large number of arrivals during her first service. Thus, such a customer might want to renege if they are allowed. But once they choose to remain, the residual time until their next service has an Erlang distribution, which has an increasing hazard rate. Thus, if it is worth remaining in the system, it is worth waiting until the next service. In the second part of this paper, we assume that customers are allowed to renege every time they join the end of the queue according to the same threshold strategy with which they choose to join the system. That is, if customers choose to join the system if and only if the number of customers in the system is less than or equal to some threshold value, then they will leave the system after a service failure  if and only if the number exceeds the same threshold value. With matrix analytic methods, we can easily compute the Nash equilibrium threshold when reneging is permitted, and compare it with the equilibrium threshold value in the non-reneging case. We show that the customers' equilibrium threshold value when reneging is allowed is greater. However, for some parameter values, their expected payoff can decrease.

The paper is organised as follows. In Section 2 we introduce the basics of the $M/M/1$ feedback queue, and precisely define our threshold joining strategy, which is specified by a real-valued threshold. We also derive an analytical expression for a tagged customer's position-dependent expected sojourn time if the other customers always choose to join. 
In Section \ref{sec: NRcase} we obtain numerically the expected sojourn time and the expected payoff of a tagged customer conditioned on her joining position and the threshold strategy used by others, via matrix analytic methods. Then we propose a threshold Nash equilibrium strategy. In Section 4 we assume that customers are allowed to leave after joining and their reneging threshold is the same as the one with which they choose to join the system. We compute the Nash equilibrium threshold when reneging is permitted and compare it with that in the non-reneging case. In Section 5 we present two paradoxes observed in the non-reneging and the reneging case. In Section 6 we analyse the optimal social welfare in both the non-reneging and the reneging cases, and prove that allowing reneging does not change the socially optimal threshold and optimal social welfare.

\section{Preliminaries} \label{sec: Pre}

\subsection{Joining strategies} \label{sec: Pre: Jpolicies}

We assume that the queue starts at time $0$ with an initial number of customers according to a distribution $\pi(0)$ which is supported on the nonnegative integers. The number of customers in the system is observable to any arriving customer before she decides to join or not to join. 
For $r = 1,2,\ldots$, let $u_r$ be a function that maps the numbers $1,2, \ldots$ to the interval $[0,1]$ such that $u_r(i)$ is the probability that the $r$th arriving customer chooses to join if there are $i-1$ customers in front of her (including the one in service), which would mean that she starts in position $i$.
We call the function $u_r$ the {\it joining strategy} for customer $r$ and $\bm{u}^\infty \equiv (u_1,u_2,\ldots )$ the {\it joining strategy profile} for the population. If $u_r(i)$ depends only on $i$, then the joining strategy is {\it symmetric} in which case, $\bm{u} = \{u,u,\cdots\}$ (see \citet[p3]{HH03}).

Next, we introduce the definition of a threshold strategy. This threshold strategy was first proposed in \citet{H96}, and was also used in \citet{HH97}.

\begin{definition} {\rm (symmetric threshold strategy)}. \label{D1}
	For any $x \in \mathbb{R}^+$, the symmetric threshold strategy with threshold value $x$ has components
	\begin{equation}
	u^{(x)}(i)\equiv
	\begin{cases}
	1  & \text{if}\ \, i \leq n  \\[+6pt]
	p  & \text{if}\ \, i=n +1 \\[+6pt]
	0  & \text{if}\ \, i \geq n +2 \,,
	\end{cases}
	\end{equation} 
\end{definition}
where $n \equiv \lfloor x \rfloor, p \equiv x-n$.
A customer who adopts threshold $x$ always chooses to join at a position which is less than or equal to $x$. She chooses to join at position $\lfloor x \rfloor+1$ with probability $x - \lfloor x \rfloor$, and refuses to join at any position greater than $\lfloor x \rfloor+1$. 
In their unpublished report \cite[Theorem 6]{BC13}, Brooms and Collins claimed that any symmetric equilibrium joining strategy must be a threshold strategy. However their proof lacks detail, so we are going to treat this result with caution. If it is correct then our threshold strategy in Theorem 1 is the unique symmetric subgame perfect equilibrium strategy.

\subsection{Basics of a single-server feedback queue} \label{sec: Pre: BasicsFeedback}

For the single-server feedback queue in Figure \ref{F1}, in the time interval $[0, \infty)$, we denote by $\xi(t)$ and $\tau_r$ the number of customers in the system at time $t$ and the arrival time of the $r$th customer, respectively. Then $\xi_r:=\xi(\tau_r)$ is the position at which the $r$th customer joins the system where, when $\xi_r = 1$, the customer immediately goes into service. 

To work out the Nash equilibrium strategy, we arbitrarily select a customer as our tagged customer, and calculate her optimal response based on different strategies adopted by others. We are interested in the symmetric Nash equilibrium strategy, that is the strategy which is the best response when others use it too.

We denote the total sojourn time of the tagged customer in the system when the other customers all use threshold $x$ by $W^{(x)}$. Consistent with this notation, $W^{(\infty)}$ is the total sojourn time of a tagged customer in the system when all the other arriving customers always join and are not allowed to renege later.

From \citet[Theorem 1]{T63}, if $\lambda < q \mu$, then when all customers always join and are not allowed to renege later, the process $\lbrace \xi(t), 0 \leq t < \infty \rbrace$ has a unique stationary distribution 
\[
\pi_i:= \lim\limits_{t \rightarrow \infty}\mathbb{P} \lbrace \xi(t) = i \rbrace =  \left(1-\frac{\lambda}{q\mu}\right)\left(\frac{\lambda}{q\mu}\right)^i  \, (i=0,1,\cdots)\,.
\]

Furthermore, \citet[Section VI]{T63} gave the Laplace-Stieltjes transform of the unconditional stationary waiting time. We use similar techniques to obtain the conditional expected sojourn time given the joining position of each customer. In the stationary regime, for $i = 1,2,...$, let
\begin{align}
&P_i(w):= \mathbb{P}\{W^{(\infty)} \leq w \, , \, \xi_r = i\} \\
&\Pi_i(s):= \int_{0}^{\infty}e^{-sw}dP_i(w) \,.
\end{align}
Then for $|z| \leq 1$, $\mathfrak{R}(s) \geq 0$, 
\begin{equation}
U(s,z) := \sum_{i=1}^{\infty} \Pi_i(s)z^i = 
\left(1-\frac{\lambda}{q \mu}\right)\sum_{k = 1}^{\infty}\frac{q (1-q)^{k-1}}{a_k(s,z)-b_k(s,z)} \,, \label{eq6}
\end{equation}
where
\begin{equation}
\begin{bmatrix}
a_k(s,z)  \\
b_k(s,z) 
\end{bmatrix}
=
\begin{bmatrix}
\frac{\displaystyle \mu+\lambda+s}{\displaystyle \mu}      &  -q \\
\frac{\ds \lambda}{\ds \mu}  & (1-q)
\end{bmatrix}^k
\begin{bmatrix}
1  \\
\frac{\ds\lambda z}{\ds q\mu}
\end{bmatrix} \,.
\end{equation}

To obtain $\displaystyle \int_{0}^{\infty}w \, dP_i(w)$, we take the derivative of $U(s,z)$ with respect  to $s$ and set $s=0$.
\begin{align}
\sum_{i=1}^{\infty} \left(\int_{0}^{\infty}w \, dP_i(w)\right) z^i
&= -\frac{\partial U(s,z)}{\partial s} \mid_{s=0}  \\
&=\left(1-\frac{\lambda}{q \mu}\right)\frac{q\mu((1-q)\lambda z+(q-2)q\mu)}{(\lambda z-q\mu)^2((q-1)\lambda-(q-2)q\mu)} \\
&=\sum_{i = 1}^{\infty} \left(1-\frac{\lambda}{q \mu}\right) \, \frac{i+1-q}{(q\mu)^i((q-1)\lambda-(q-2)q\mu)}z^i \,. \label{eq3}
\end{align}
Hence, the stationary expected sojourn time of a tagged customer if she joins at position $i$, and all other customers always choose to join upon arrival is
\begin{equation}
w_{i,i}^{(\infty)}:=\mathbb{E}\left(W^{(\infty)} \mid \xi_r = i\right) = \frac{\ds\int_{0}^{\infty}w \, dP_i(w)}{\ds \pi_i} = \frac{i+1-q}{((q-1)\lambda-(q-2)q\mu)}  \,. \label{closed_form_waiting_time}
\end{equation}

\section{The Case When Customers Cannot Renege}  \label{sec: NRcase}
\subsection{Expected payoff} \label{sec: NRcase: EP}
In this paper, we assume that customers are homogeneous which means they value receiving service identically and they place the same per unit time value on their waiting, and we focus on symmetric threshold strategies defined in Definition \ref{D1}. When a customer arrives and sees $j-1$ customers already in the system, she will join the queue at the $j$th place. When every customer adopts threshold $x$ and the system starts with less than $\lceil x \rceil+1$ customers, the tagged customer, upon arrival, can observe at most $\lceil x \rceil$ people in the system. If she chooses to join, her position is at most $\lceil x \rceil + 1$. 

Let $w^{(x)}_{i,j}$ be the expected remaining time until the tagged customer departs the system, if there are $j$ customers in the system, she is in position $i\leq j$ and all the other customers use threshold $x$. So if a customer joins in position $j$, her expected sojourn time will be $w_{j,j}^{(x)}$. On the other hand, when she leaves the queue she will obtain a reward $R_0$ and her expected payoff when she is in position $i$ there are $j$ customers in total and other customers are using threshold $x$ is thus $z_{i,j}^{(x)} \equiv E(R) - w_{i,j}^{(x)} =  R_0 - w_{i,j}^{(x)}$.

We shall show that the vector 
\[
\bm{w}^{(x)} = \left(w_{1,1}^{(x)}, w_{1,2}^{(x)}, w_{2,2}^{(x)}, \ldots, w_{1,\lceil x \rceil +1}^{(x)}, \ldots, w_{\lceil x \rceil,\lceil x \rceil +1}^{(x)}, w_{\lceil x \rceil +1,\lceil x \rceil +1}^{(x)}\right)^T  
\] 
satisfies a version of Poisson's equation. In Section \ref{sec: Rcase} where we consider a model with reneging, customers do not always get the reward, and we proceed by writing Poisson's equation for the expected payoff $z_{i,j}^{(x)}$ directly.

To compute $\bm{w}^{(x)} $, we construct a continuous-time quasi-birth-and-death process (QBD) on the state space $\mathcal{S} = \left\{(i,j) \, : \, 1 \leq i \leq j \leq \lceil x \rceil + 1 \right\}$ with its level $j$ denoting the total number of customers including the customer in service in the system, and its phase $i$ denoting the position of the tagged customer. Then we construct the embedded discrete-time QBD obtained by observing this continuous-time Markov chain at its transition points and write $w_{i,j}^{(x)}$ conditioning on the first transition out of state $(i,j)$ in \eqref{eq:wij}. Specifically, the expected time until the next transition is $
\frac{1}{\lambda+\mu}$. The next transition is an arrival with probability $ \frac{\lambda}{\lambda+\mu}$. When $j < \lfloor x \rfloor$, the arriving customer joins the system with probability $1$; when $j = \lfloor x \rfloor$, the arriving customer joins the system with probability $p$; when $j = \lfloor x \rfloor+1$ or $\lfloor x \rfloor+2$, the arriving customer balks.

The next transition is a service completion with probability $\frac{\mu}{\lambda+\mu}$, after which a customer leaves the system with probability $q$ and joins the end of the system with probability $1-q$. Hence, if the customer in service is the tagged one ($i=1$), when she finishes her service, her future sojourn time is $0$ with probability $q$, otherwise, her next position is $j$. When the customer in service is not the tagged one, the position of the tagged customer decreases by $1$, the total number of customers decreases by $1$ with probability $q$ but stays unchanged with probability $1-q$. From the aforementioned reasoning,
\begin{align}  \label{eq:wij}
&\lefteqn{w_{i,j}^{(x)} \, = }  \\
& \frac{1}{\lambda+\mu}  + \frac{\lambda}{\lambda+\mu} \,  \left( w_{i,j+1}^{(x)}  \, \mathbbm{ 1 }_{\lbrace j <  \lfloor x \rfloor \rbrace}  + \left( p \, w_{i,j+1}^{(x)} + (1-p) \, w_{i,j}^{(x)}  \right) \, \mathbbm{ 1 }_{\lbrace  j =  \lfloor x \rfloor \rbrace} + w_{i,j}^{(x)} \,  \mathbbm{ 1 }_{\lbrace  j = \lfloor x \rfloor+1, \lfloor x \rfloor+2 \rbrace} \right) \, \notag \\
&\, + \, \frac{\mu}{\lambda+\mu}\, \left( (1-q) \, w_{j,j}^{(x)} \, \mathbbm{ 1 }_{\lbrace  i = 1 \rbrace} + \left( q \, w^{(x)}_{i-1,j-1} + (1-q) \, w_{i-1,j}^{(x)} \right) \, \mathbbm{ 1 }_{\lbrace  i > 1 \rbrace} \right) \,.\notag
\end{align}
Thus, we can obtain $\bm{w}^{(x)}$ by solving Poisson's equation  
\begin{equation} 
\left(I-P^{(x)} \right)
\bm{w}^{(x)}  \, 
=  \, \frac{1}{\lambda+\mu} \, \bm{e}, \label{eq: poisson1}
\end{equation}
where $P^{(x)}$ is defined in Appendix \ref{appendix:1.1}, and $\bm{e}$ denotes a vector of $1$'s of the appropriate size.

We have shown that $w_{i,j}^{(x)}$ can be obtained by solving a system of linear equations. However, the number of equations is quadratic in $\lfloor x \rfloor$. Thus, it is necessary to come up with an efficient way of carrying out the calculation. Equation (\ref{eq: poisson1}) is Poisson’s equation for a level dependent QBD, where the defining matrices $A^{(j)}_1, A^{(j)}_0, A^{(j)}_{-1}$ are given in Appendix \ref{appendix:1.1}. Due to the special structure of QBDs, we propose Algorithm \ref{alg_poisson} to solve $\bm{w}^{(x)}$ based on the methodology in \citet{DLL13}. See \citet[Chapter 12]{LR99} for a detailed explanation of the matrices $\Gamma^{(j)}, U^{(j)}$ and $G^{(j)}$ for a level dependent QBD which are used in Algorithm \ref{alg_poisson}, noting that the matrix $\Gamma^{(j)}$ in this paper has the same meaning as matrix $R^{(j)}$ in \cite{LR99}. We use $\Gamma$ to differentiate it from the reward $R$ that is obtained by the customers after they leave the service.
\begin{algorithm} 
	\caption{}\label{Poisson}
	\begin{algorithmic}[1]
		\Procedure{Calculate $U^{(j)}$, $\Gamma^{(j)}$, $G^{(j)}$}{} \Comment{The $U^{(j)}, \Gamma^{(j)}, G^{(j)}$ of $P^{(x)}$}
		\State $U^{(\lceil x \rceil+1)} \gets A_0^{(\lceil x \rceil+1)}$
		\State $\Gamma^{(\lceil x \rceil+1)} \gets A_1^{(\lceil x \rceil)} \, (\mathbf{I}- U^{(\lceil x \rceil+1)})^{-1}$
		\State $G^{(\lceil x \rceil+1)} \gets (\mathbf{I}- U^{(\lceil x \rceil+1)})^{-1} \, A_{-1}^{(\lceil x \rceil+1)}$
		\For {$j = \lceil x \rceil:2$}
		\State $U^{(j)} \gets A_0^{(j)} + A_1^{(j)}\, G^{(j+1)}$
		\State $\Gamma^{(j)} \gets A_1^{(j-1)} \, (\mathbf{I}- U^{(j)})^{-1}$
		\State $G^{(j)} \gets (\mathbf{I}- U^{(j)})^{-1} \, A_{-1}^{(j)}$
		\EndFor
		\State \textbf{end}
		\EndProcedure
		\Procedure{Poisson's Equation}{$U^{(j)}, \Gamma^{(j)}, G^{(j)}, \frac{1}{\lambda+\mu}$}
		\State $y(1) \gets 0$
		\For {$j = 2:\lceil x \rceil+1$}
		\State $y(j) \gets \frac{1}{\lambda+\mu}(\mathbf{I} - U^{(j)})^{-1} \, (\bm{e}_{j}+ \sum_{k = j}^{\lceil x \rceil} \Pi_{l=j+1: k+1} \, \Gamma^{(l)} \, \bm{e}_{k+2}) + \, G^{(j)} \, y(j-1)$
		\EndFor
		\State \textbf{end}
		\State $y(1) \gets \frac{1}{\lambda + \mu} + A_1^{(1)}\, y(2)$
		\State $w(1) = \frac{y(1)}{1-(A_0^{(1)} \, + \, A_1^{(1)} \, G^{(2)})}$ \Comment{Expected sojourn time}
		\For {$j = 2: \lceil x \rceil+1$}
		\State $w(\frac{j(j-1)}{2}+1: \frac{j(j-1)}{2}+j) = y(j) \, + \, \Pi_{l=j:2} \, G^{(l)} \, w(1)$
		\EndFor
		\State \textbf{end}
		\State \textbf{return} $w$
		\EndProcedure
	\end{algorithmic} \label{alg_poisson}
\end{algorithm}

\begin{figure}
	\centering
	\includegraphics[width = 0.7\textwidth]{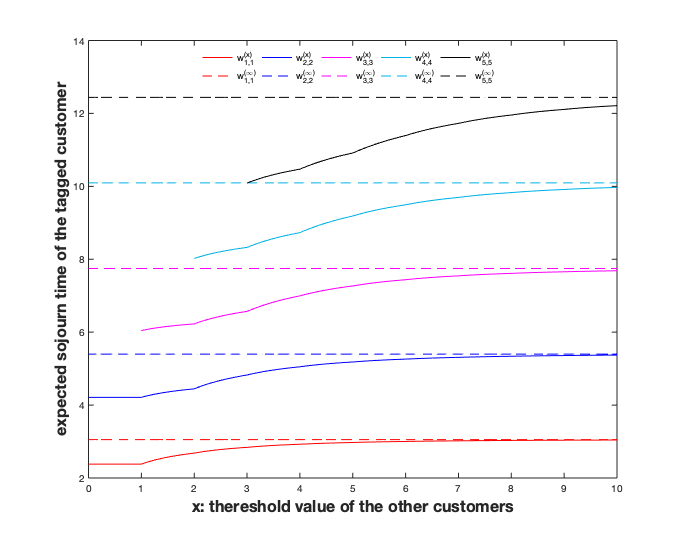}
	\caption{Expected sojourn time of the tagged customer ($\lambda = 0.4, \mu = 0.6, q = 0.7$).} \centering\label{WT}
\end{figure}

We plot $w_{j,j}^{(x)}$ for $0 \leq x \leq 10$ in Figure \ref{WT}. Several observations can be made. First, $w_{j,j}^{(x)}$ exists only when $\lceil x \rceil \geq j-1$. The reason follows from the explanation at the beginning of this section that the tagged customer cannot be in a position greater than $\lceil x \rceil + 1$. Second, $w_{j,j}^{(x)}$ increases in $j$ for any $1 \leq j \leq \lceil x \rceil+1$, and increases in $x$ when $x \geq 1$. This property was proved in \citet{BC13} via coupling, and their proof works for $GI/G/1$ feedback queues. 
For an $M/M/1$ feedback queue, we propose an alternative proof in Lemmas \ref{lemmaWI} and \ref{lemmaWx}. Third, when $x \in [0,1)$, as long as the tagged customer is in the system, no newly arriving customer will join the system, hence the expected sojourn time of the tagged customer is independent of $x$. Actually, from (\ref{eq:wij}), we explicitly have 
\begin{equation}
w_{1,1}^{(x)}  = \frac{1}{\mu q} \qquad w_{2,2}^{(x)} = \frac{3-q}{\mu q (2-q)} > \frac{1}{\mu q} \qquad 0 \leq x \leq 1 \,.
\end{equation}
Finally, as expected, $w_{j,j} ^{(x)}$ approaches $w_{j,j} ^{(\infty)}$ as $x$ increases. Our results are stated in Lemmas 1 and 2 below, the proofs of which appear in Appendix \ref{appendix:proof of lemma1 and 2}.

\begin{lemma} \label{lemmaWI}
	$w^{(x)}_{j,j}$ is increasing in $j$ for $1 \leq j \leq \lceil x \rceil +1 \,.$
\end{lemma}

\textbf{Remark} At the expense of making the calculation more intricate, we can prove that $w_{j,j}^{(x)}$ is strictly increasing in $j$. We omit the details.

\begin{lemma} \label{lemmaWx}
	For any two threshold policies $x_1$ and $x_2$ with $x_1 < x_2$, 
	\begin{equation}
	w_{i,j}^{(x_1)} <  w_{i,j}^{(x_2)} \qquad 1 \leq i \leq \lceil x_1 \rceil +1 \,.
	\end{equation}
\end{lemma}

\subsection{The Nash equilibrium threshold} \label{sec: NRcase: NashE}
In Lemmas \ref{lemmaWI} and \ref{lemmaWx}, we have proved that $w_{j,j}^{(x)}$ is increasing in $j$ and $x$, so $z_{j,j}^{(x)}$ is decreasing in $j$ and $x$. We know from the beginning of Section \ref{sec: NRcase} that the position where the tagged customer can join is at most $\lceil x \rceil + 1$ if the other customers use threshold $x$ and the system starts with less than $\lceil x \rceil$ customers. If we refer to the highest position that the tagged customer is willing to join, when others use threshold $x$, as the best response, and let $\mathcal{BR}(x)$ denote it, then $\mathcal{BR}(x) = \max\{ j: z^{(x)}_{j,j} \geq 0, \, 1 \leq j \leq \lceil x \rceil +1  \}$.

If $R_0$ is big, then there will be values of $x$ for which  $z^{(x)}_{\lceil x \rceil +1, \lceil x \rceil +1} \geq 0$ and so $\mathcal{BR}(x) = \lceil x  \rceil + 1$. On this part of the domain, $\mathcal{BR}(x)$ is (obviously) an increasing step function. However, as $x$ increases, there must be a value $x^*$ for which $z^{(x^*)}_{\lceil x^* \rceil +1, \lceil x^* \rceil +1} < 0$. To see this, observe that a customer arriving to position $j$ must wait for at least $j$ services and so $w^{(x)}_{\lceil x \rceil +1, \lceil x \rceil +1} \geq (\lceil x \rceil +1)/\mu$ and so when $x > R_0\mu - 1$,
\begin{eqnarray}
z^{(x)}_{\lceil x \rceil +1, \lceil x \rceil +1} & = & R_0 - w^{(x)}_{\lceil x \rceil +1, \lceil x \rceil + 1} \\
         & \leq & R_0 - (\lceil x \rceil +1)/\mu \\
  & < & 0.
\end{eqnarray}
For $x>x^*$, $BR(x) < \lceil x \rceil + 1$ and, on this part of the domain Lemma 2 ensures that $BR(x)$ is a monotone decreasing step function.

There are now two possibilities 
\begin{itemize}
\item there is an integer $m$ such that $BR(m) = m$, or
\item there is an integer $m$ such that $BR(m) = m+1$ and $BR(m+1) \leq m$,
\end{itemize}
These are illustrated in Figure 3.2(a) and in Figures 3.2(b) and (c) respectively.

\begin{figure}[h]
	\subcaptionbox{$\mathcal{BR}(m) = m$}%
	{\includegraphics[width=0.29\linewidth]{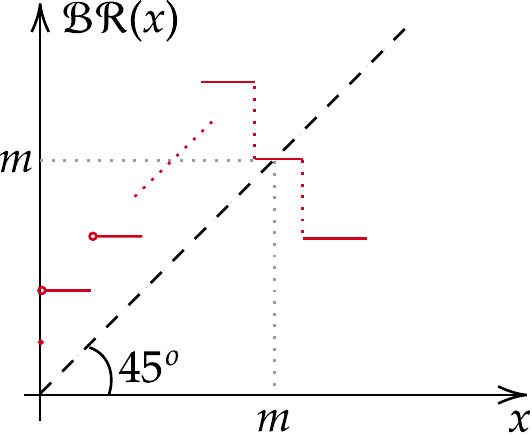}}
	\subcaptionbox{$\mathcal{BR}(m) = m+1,\, \mathcal{BR}(m+1) = m$}%
	{\includegraphics[width=0.34\linewidth]{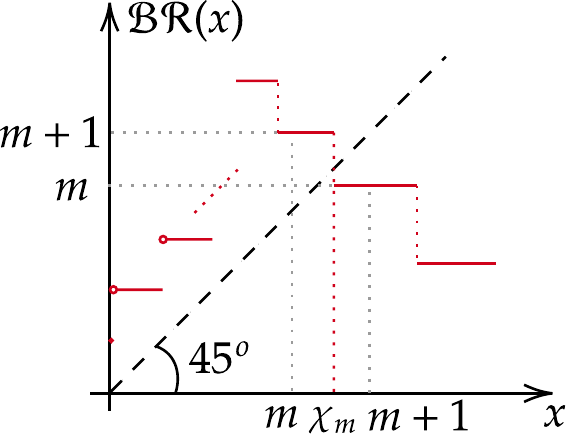}}	
	\hspace{\fill}
	\subcaptionbox{$\mathcal{BR}(m) = m+1,\, \mathcal{BR}(m+1) < m$}%
	{\includegraphics[width=0.34\linewidth]{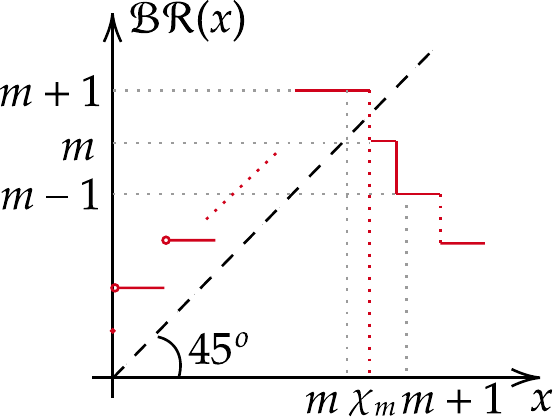}}	
	\caption{Best Response.}
	\label{fig:BR}
\end{figure}
 
\begin{figure}
	\centering
	\includegraphics[width = .7\textwidth]{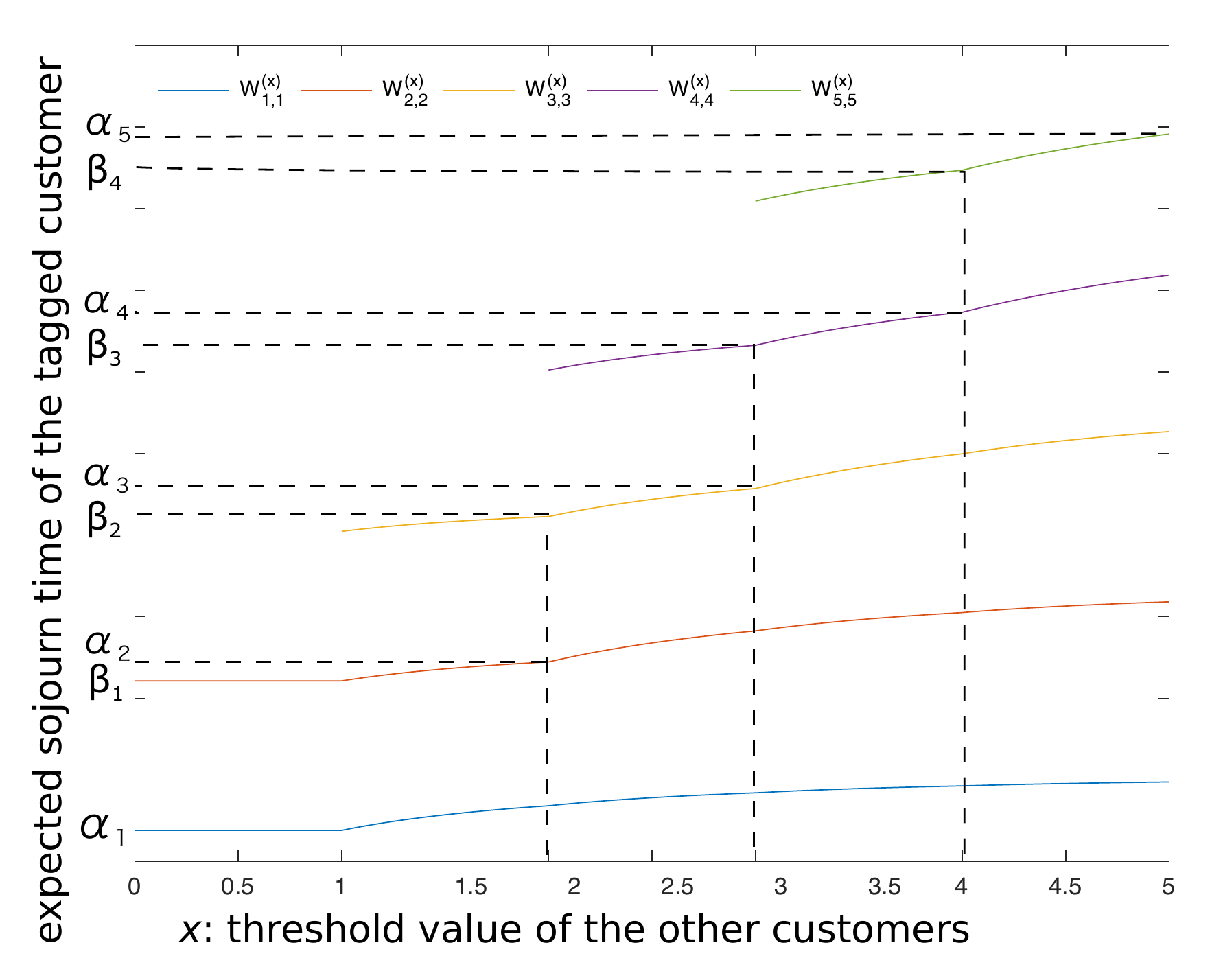}
	\caption{} \centering \label{NEp}
\end{figure}

For the purpose of presenting the Nash equilibrium, for $m = 1, 2, \cdots$, let $\alpha_m = w_{m,m}^{(m)}, \beta_m = w_{m+1,m+1}^{(m)}$ (see Figure \ref{NEp}). Let $\chi_m$ be the solution to $w_{m+1,m+1}^{(\chi_m)} = R_0$. We prove in the following lemma that $\chi_m$ exists and is unique.
 \begin{lemma}
 	There exists a unique $\chi_m$ such that $w_{m+1,m+1}^{(\chi_m)} = R_0$.
 \end{lemma}
Proof. 
For $x \in (m,m+1)$, from Equation \eqref{eq: poisson1}, 
\begin{equation}
	w_{m+1,m+1}^{(x)} = \frac{1}{\lambda+\mu} \left( \left( I-P^{(x)} \right)^{-1}\, \bm{e} \right)_{\frac{(m+1)(m+2)}{2}} \,,
\end{equation}
where $P^{(x)}$ is defined in Appendix \ref{appendix:1.1}. The matrix $P^{(x)}$ is substochastic for any $x$, as the sum of the $\left(\frac{j(j+1)}{2}+1\right)$th row of $P^{(x)}$ is less than $1$ for $j = 1, 2 \, \ldots, \lceil x \rceil + 1$. From a Corollary to the Perron-Frobenius Theorem (\citet[page 8]{S06}), $|r^{(x)}| < 1 $ for any eigenvalue $r^{(x)}$. Thus, any real eigenvalue $1-r^{(x)}$ of $I-P^{(x)}$ must be greater than 0. Hence $|I-P^{(x)}| \ne 0$. 

Next, we write $x$ as $m+p$. From its expression, $P^{(x)}$ is continuous in $p$, so is $I-P^{(x)}$. Since the entries of the inverse matrix can be written as rational functions of the entries of the original matrix, and the denominators of these rational functions are non-zero for all $x$,  $w_{m+1,m+1}^{(x)}$ is continuous in $p$. Hence, $w_{m+1,m+1}^{(x)}$ is continuous in $x$ for $x \in (m,m+1)$. Also, it follows from Lemma \ref{lemmaWI} and \ref{lemmaWx} that $\beta_m > \alpha_m$ and $\alpha_{m+1} > \beta_m$, respectively. If $\beta_m < R_0 < \alpha_{m+1}$, then
\[
\beta_m  = w_{m+1,m+1}^{(m)} < R_0 < w_{m+1,m+1}^{(m+1)} = \alpha_{m+1}  \,.
\]
Due to the fact that $w_{m+1,m+1}^{(x)}$ continuous and strictly increasing in $x \in (m,m+1)$, there is a unique $\chi_m \in (m,m+1)$ such that $w_{m+1,m+1}^{(\chi_m)} = R_0$.

 We describe the Nash equilibrium strategy for the feedback queueing system in the following theorem.
\begin{theorem}\label{Theorem:NE}
	There exists an equilibrium threshold strategy with threshold value 
	\begin{equation}\label{NE}
	x_e=
	\begin{cases}
	0 & \text{if} \ R_0 < \alpha_1 \,,\\[+6pt]
	\zeta_0 & \text{if} \ R_0 = \alpha_1 \,,\\[+6pt]
	m &  \text{if}	\ \alpha_m	\leq R_0  \leq	\beta_m \quad m = 1,2, \ldots,\\[+6pt]
	\chi_m & \text{if}	\ \beta_m 	<  R_0 	<	\alpha_{m+1} \quad m = 1,2,\ldots \,,
	\end{cases}
	\end{equation}
	where $\zeta_0 \in [0,1]$.
\end{theorem}

Proof. A customer will choose to join the queue if and only if her reward can fully bear her expected sojourn cost.
\begin{itemize}
	\item When $R_0< \alpha_1$, even if the tagged customer is the only one in the system, her expected sojourn time $\alpha_1 > R_0$. Thus, her best option is balking. The same analysis works for the other customers. So balking is the Nash equilibrium strategy.
	\item When $R_0= \alpha_1$, if the other customers are all using threshold $\chi_0 \in \left[0,1\right]$, there is at most one customer in the system when the tagged customer arrives. For the tagged customer, when she observes one person in the system, her best response is balking as her expected sojourn time $\beta_1 > \alpha_1 = R_0$. When she observes that the system is empty, her expected payoff is zero. Thus, she is indifferent between joining an empty system and balking. Actually, she gains nothing by either strategy. The same analysis works for any other customer, so any threshold strategy with threshold value $\chi_0 \in \left[0,1\right]$ is a Nash equilibrium strategy.  
	\item When $\alpha_m \leq R_0 \leq \beta_m$, the tagged customer's expected sojourn time satisfies $w_{m,m}^{(m)} \leq R_0\leq w_{m+1,m+1}^{(m)}$, so                                                                                                                                                        
	\begin{equation}             
	z_{m+1,m+1}^{(m)} \leq 0 \leq z_{m,m}^{(m)}\,.                                                           
	\end{equation}                                                                                                                                  
	Hence the tagged customer's best response is $m$ when others adopt threshold $m$. So threshold $m$ is a Nash equilibrium strategy.                                                                                                                                                                                                           
	\item When $\beta_m	<  R_0 < \alpha_{m+1}$,                                                                                                                                                                                                                                                                                                                                                                                              if other customers all adopt threshold $\chi_m$, the tagged customer gains nothing when she joins at $m+1$, so her best response is any threshold strategy with threshold value between $m$ and $m+1$ (including $\chi_m$). Thus, $\chi_m$ is the Nash equilibrium threshold. \hfill $\square$
\end{itemize}                                                                                                                                                                                                                                                                                                                 From Theorem \ref{Theorem:NE}, $m$ is either the Nash equilibrium threshold or the integer part of it. The Nash equilibrium threshold is not an integer when $R_0 \in (\beta_m, \alpha_{m+1})$. Figure \ref{fig:BR}(b) represents this case with $\mathcal{BR}(m) = m+1$ and $\mathcal{BR}(m+1) = m$. Figure \ref{fig:BR}(c) depicts the case with $\mathcal{BR}(m) = m+1$ and $\mathcal{BR}(m+1) = m-1$. In both cases, the tagged customer is indifferent between $m$ and $m+1$ when others use a threshold between $m$ and $m+1$, and the conclusion of Theorem 3.1 holds.

Intuitively speaking, a Nash equilibrium is said to be {\em  evolutionarily stable} if it cannot be invaded by any alternative strategy that is initially rare (see \citet{S86}).                                                                                                                         
\begin{definition}  
	\textbf{Evolutionarily stable strategy (ESS).} A Nash equilibrium strategy $x$ is said to be an ESS if either (i) $x$ is the unique best response against itself or (ii)  for any $x' \ne x$ which is a best response against $x$, $x$ is better than $x'$ as a response to $x'$ itself. That is, with $U(x',x)$ denoting a customer's expected payoff when she uses $x'$ and others use $x$, for all $x'\not = x$, either
	\begin{align}                                                  
	&U(x, x) > U(x', x)\,, \, \mbox{or} \\                                                                                                                                                                                                                                                                                                       
	&U(x, x) = U(x', x) \quad \mbox{and} \quad U(x, x') > U(x', x') \,.                                                                                                                                                                                                                                                                          
	\end{align}                                                                                                                       
\end{definition}                                                                                                                                                                                                                                                                                                                    
To show that the Nash equilibrium startegy with threshold value $x_e$ is an ESS, we first define the total expected payoff of a tagged customer who adopts threshold $x$ when the other customers all adopt threshold $x'$.          
\begin{equation} 
\label{eq:TotalEP}
U(x,x'): = \sum_{i=1}^{\lfloor x \rfloor \wedge (\lceil x' \rceil + 1)} \pi^{(x')}_{i-1} \, z_{i,i}^{(x')} +  (x - \lfloor x \rfloor) \, \pi^{(x')}_{\lfloor x \rfloor} \, z_{\lfloor x \rfloor+1,\lfloor x \rfloor+1}^{(x')} \, \mathbbm{1}_{\{\lfloor x \rfloor < \lceil x' \rceil + 1\}}  \,,
\end{equation}
where $\pi^{(x)}_j \,\, 0 \leq j \leq \lceil x \rceil$, is the stationary distribution of the number of customers in the system where everyone adopts threshold $x$. We prove the $x_e$-threshold strategy is an ESS in the following corollary.
\begin{corollary}
	The threshold strategy with threshold value $x_e$ is an ESS when $R \neq \alpha_1$.
\end{corollary}
We have already proved that when other customers adopt threshold strategy $x_e$, there is no better strategy than $x_e$ for the tagged customer, that is $U(x_e, x_e) \geq U(x, x_e)$. 
\begin{itemize}
	\item When $R_0< \alpha_1$,  $U(0,0) = 0 > U(x,0)$ for any $x>0$. Thus, balking is an ESS.
	\item When $R_0= \alpha_1$,  for any $0 \leq \chi_0,\chi_0'\leq 1$, $U(\chi_0,\chi_0) = U(\chi_0',\chi_0) = 0$ and $U(\chi_0,\chi_0') = U(\chi_0',\chi_0') = 0$. Thus, $\chi_0 \in [0,1]$ is not an ESS.
	\item When $\alpha_m \leq R_0 \leq \beta_m$, $U(m,m) > U(m',m)$ for any $m' \ne m$. Thus, $m$ is an ESS.
	\item When $\beta_m	<  R_0	< \alpha_{m+1}$, it follows from the definition of $\chi_m$ that $z_{m+1,m+1}^{(\chi_m)} = R_0 - w_{m+1,m+1}^{(\chi_m)} = 0$, so for any $\chi_m' \in [m, m+1]$
	\begin{equation}
	U(\chi_m', \chi_m) = U(\chi_m, \chi_m) = \sum_{i=1}^{m} \pi^{(\chi_m)}_{i-1} \, z_{i,i}^{(\chi_m)}  \,.
	\end{equation}
	Furthermore, when $m \leq \chi_m' <\chi_m < m+1$, it follows from the fact that $z_{m+1,m+1}^{(x)}$ is decreasing in $x$ that
	\begin{equation}
	z_{m+1,m+1}^{(\chi_m')} >0 \,.
	\end{equation}
	Since $\lfloor \chi_m' \rfloor = \lfloor \chi_m \rfloor = m$, the first summations in $U(\chi_m', \chi_m')$ and $U(\chi_m, \chi_m')$ are equal. However, $(\chi_m' - \lfloor \chi_m' \rfloor) < (\chi_m - \lfloor \chi_m \rfloor)$, hence the second term in $U(\chi_m',\chi_m') $ is less than the second term in $U(\chi_m, \chi_m')$. So $U(\chi_m',\chi_m') < U(\chi_m,\chi_m')$. Similarly, when $\chi_m < \chi_m' \leq m+1$, $z_{m+1,m+1}^{(\chi_m')} <0$ but $(\chi_m - \lfloor \chi_m \rfloor) < (\chi_m' - \lfloor \chi_m' \rfloor)$. Following similar lines, we have $U(\chi_m',\chi_m') < U(\chi_m,\chi_m')$. Thus, $\chi_m$ is an ESS.
$ \hfill \square $
\end{itemize}
\begin{figure} 
	\centering
	\includegraphics[width=10cm]{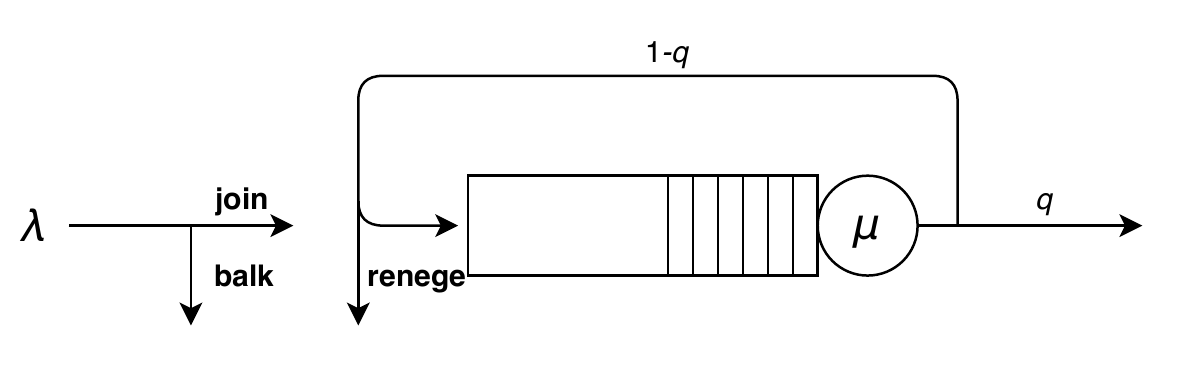}
	\caption{An M/M/1 feedback queue with strategic customers when reneging is allowed.} \label{RC}
\end{figure}

\section{The Case When Customers Can Renege} \label{sec: Rcase}
Every time a customer rejoins at the end of the queue due to a service failure, it is possible that her conditions have deteriorated with time. Hence customers might have an incentive to {\it renege}, that is depart from the queue, when their service fails. Figure \ref{RC} is an illustration of an $M/M/1$ feedback queue when reneging is allowed. In this section, we focus on the Nash equilibrium threshold when the customers are allowed to renege, and compare it with the equilibrium threshold when they cannot renege. In order to make comparisons between the two cases, we abbreviate the non-reneging case as the $N$-case and the reneging case as the $R$-case.
\subsection{The expected payoff} \label{sec: Rcase:EP}

In our model, every time a customer rejoins the end of the queue, she faces a similar situation as that when she first chooses to join. Thus, we restrict our attention to policies where the customer must use the same threshold when she chooses to balk or renege.

In contrast to Section \ref{sec: NRcase}, the expected payoff of the tagged customer is affected by her future reneging decisions. In particular, if she chooses to renege, she will not receive the reward $R_0$. So we use $\hat{z}^{(x_{tag}, x)}_{i,j}$ to denote the tagged customer's expected payoff, which is the difference between the expected reward and her expected sojourn cost, given that she is at position $i$ and uses threshold strategy $x_{tag}$, there are $j$ customers in the system, and the other customers all adopt threshold $x$. It will turn out that the relevant value of $x_{tag}$ that we need to consider for the purpose of calculating the Nash equilibrium occurs when $x _{tag} = \lfloor x \rfloor +1$. This satisfies the equation
\begin{align} 
& \hat{z}_{i,j}^{(\lfloor x \rfloor +1,x)}  = -\frac{1}{\lambda+\mu}   \\
&+ \frac{\lambda}{\lambda+\mu} \, \left( \hat{z}_{i,j+1}^{(\lfloor x \rfloor +1,x)} \, \mathbbm{1}_{\lbrace  j < \lfloor x \rfloor \rbrace} +\left(p \, \hat{z}_{i,j+1}^{(\lfloor x \rfloor +1,x)} + (1-p) \, \hat{z}_{i,j}^{(\lfloor x \rfloor +1,x)} \right) \, \mathbbm{1}_{\lbrace j = \lfloor x \rfloor \rbrace} \right. \notag \\
& \left. + \hat{z}_{i,j}^{(\lfloor x \rfloor +1,x)} \, \mathbbm{1}_{\lbrace  j = \lfloor x \rfloor +1  \rbrace} \right) \notag + \frac{\mu}{\lambda+\mu} \left( \left(q \, R_0 + (1-q) \, \hat{z}_{j,j}^{(\lfloor x \rfloor +1,x)} \right) \, \mathbbm{1}_{\lbrace i = 1  \rbrace} \right. \notag \\
& \left. + \left((q + (1-q) \, (1-p) \, \mathbbm{1}_{\lbrace j=\lfloor x \rfloor+1 \rbrace}) \, \hat{z}_{i-1,j-1}^{(\lfloor x \rfloor +1,x)} + (1-q) \, (1-(1-p) \mathbbm{1}_{\lbrace j=\lfloor x \rfloor+1 \rbrace}) \, \hat{z}_{i-1,j}^{(\lfloor x \rfloor +1,x)} \right) \right) \notag \,,
\end{align}
where $p$ is the fractional part of $x$ as defined in Definition 
\ref{D1}. Hence, we can calculate $\hat{z}_{i,j}^{(\lfloor x \rfloor +1,x)}$ via Poisson's equation
\begin{equation} \label{eq: NER1}
\left(I-\hat{P}^{(\lfloor x \rfloor +1,x)} \right) \, 
\bm{\hat{z}}^{(\lfloor x \rfloor +1,x)} \, 
=  \, \bm{{g}} \,,
\end{equation}
where the matrix $\hat{P}^{(\lfloor x \rfloor +1,x)}$ and the vector $\bm{g}$ are defined in Appendix \ref{appendix: 2.1}, and 
\begin{equation}
\bm{\hat{z}}^{(\lfloor x \rfloor +1,x)} = (z_{1,1}^{(\lfloor x \rfloor +1,x)}, z_{1,2}^{(\lfloor x \rfloor +1,x)}, z_{2,2}^{(\lfloor x \rfloor +1,x)}, \cdots, z_{\lfloor x \rfloor,\lfloor x \rfloor+1}^{(\lfloor x \rfloor +1,x)}, z_{\lfloor x \rfloor+1,\lfloor x \rfloor+1}^{(\lfloor x \rfloor +1,x)})^T \,.
\end{equation}

In Section \ref{sec: NRcase: NashE}, we derived the Nash equilibrium threshold value by finding the $m$ that satisfies the case in Figure \ref{fig:BR}. In the $N$-case, this means only $z^{(\lfloor x \rfloor)}_{\lfloor x \rfloor, \lfloor x \rfloor}$, $z^{(\lfloor x \rfloor)}_{\lfloor x \rfloor+1, \lfloor x \rfloor+1}$ and $z^{(x)}_{\lfloor x \rfloor+1, \lfloor x \rfloor+1}$ matter in calculating the Nash equilibrium, although the tagged customer can join at position $\lceil x \rceil +1$. Similarly, in the $R$-case we only care about $\hat{z}^{(\lfloor x \rfloor, \lfloor x \rfloor)}_{\lfloor x \rfloor, \lfloor x \rfloor}$, $\hat{z}^{(\lfloor x \rfloor+1, \lfloor x \rfloor)}_{\lfloor x \rfloor+1, \lfloor x \rfloor+1}$, and $\hat{z}^{(\lfloor x \rfloor+1,, x)}_{\lfloor x \rfloor+1, \lfloor x \rfloor+1}$, so we calculate $\hat{z}^{(\lfloor x \rfloor+1, \lfloor x \rfloor)}_{i,j}$ only for $1 \leq i \leq  j \leq \lfloor x \rfloor+1$.

In the $R$-case, when others use threshold $x$ and the tagged customer uses threshold $x_{tag} \geq \lceil x \rceil$, the queue size is never greater than $\lceil x \rceil$ at a time point where the tagged customer's service has failed, so the tagged customer will never renege after joining if she uses $x_{tag}$ even though other customers may do so. Hence the calculation of $\bm{\hat{z}}^{(x_{tag}, x)}$ when $x_{tag} \geq \lceil x \rceil$ can be transfered to the calculation of the expected sojourn time. If we define $\hat{w}^{(x_{tag}, x)}_{i,j}$ as the expected sojourn time of the tagged customer in the $R$-case, given that she is at position $i$ and uses threshold strategy $x_{tag}$, there are $j$ customers in the system, and the other customers all adopt threshold $x$, then when $x_{tag} = \lfloor x \rfloor +1$,
{\begin{align}
\bm{\hat{w}}^{(\lfloor x \rfloor +1, x)}=\left( w_{1,1}^{(\lfloor x \rfloor +1,x)}, w_{1,2}^{(\lfloor x \rfloor +1,x)}, w_{2,2}^{(\lfloor x \rfloor +1,x)}, 
\cdots, w_{\lfloor x \rfloor,\lfloor x \rfloor+1}^{(\lfloor x \rfloor +1,x)}, w_{\lfloor x \rfloor+1,\lfloor x \rfloor+1}^{(\lfloor x \rfloor +1,x)}\right)^T \,,
\end{align}}
satisfies a version of Poisson's equation similar to Equation \eqref{eq: poisson1}
\begin{equation} \label{eq: poisson2}
\left(I-\hat{P}^{(\lfloor x \rfloor+1, x)}\right) \, \bm{\hat{w}} = \frac{1}{\lambda+\mu} \, \bm{e} \,,
\end{equation}
and $\bm{\hat{z}}^{(\lfloor x \rfloor +1, x)} = R_0-\bm{\hat{w}}^{(\lfloor x \rfloor +1, x)}$. Similar to the $N$-case, an equilibrium strategy exists and can be computed using algorithm \ref{alg_poisson}.

We compare ${\hat{z}}^{(\lfloor x \rfloor +1,x)}_{j,j}$ and ${z}^{(x)}_{j,j}$ for $j = 1, \ldots, \lfloor x \rfloor+1$ in Lemma \ref{lemma:CompareRNR1}, the proof of which appears in Appendix \ref{appendix:lemma3}.
\begin{lemma} \label{lemma:CompareRNR1}
	When $\lfloor x \rfloor < x$, 
	\begin{align} \label{eq:compare1}
	\hat{z}^{(\lfloor x \rfloor  +1,x)}_{j,j} \geq z^{(x)}_{j,j} \quad for \quad j = 1, \ldots, \lfloor x \rfloor+1 \,.
	\end{align}
	When $\lfloor x \rfloor = x$, 
	\begin{align}  \label{eq:compare2}
	\hat{z}^{(\lfloor x \rfloor, \,x )}_{j, j} = \hat{z}^{(\lfloor x \rfloor  +1,x)}_{j,j} = z^{(x)}_{j,j} \quad for \quad j = 1, \ldots, \lfloor x \rfloor \,. 
	\end{align}
\end{lemma}

One interpretation of Lemma \ref{lemma:CompareRNR1} is as follows. When other customers adopt the threshold $x$, for a customer who never reneges, her expected payoff is higher if the other customers are allowed to renege. When $x = \lfloor x \rfloor > 0$, the number of customers in the system never exceeds $\lfloor x \rfloor$ if the tagged customer joins at a position less than $\lfloor x \rfloor+1$; if the tagged customer joins at $\left(\lfloor x \rfloor+1\right)$th position, the customer who is in service when she joins will leave the system with probability $1$: either the service will complete successfully or the customer will renege when the service fails. 
Thus if the tagged customer joins at position $\lfloor x \rfloor+1$, then she is better off when others can renege, but there is no difference between the $N$-case and the $R$-case when the position at which the tagged customer joins is less than $\lfloor x \rfloor+1$.

\textbf{Remark}. 
We can prove the strict inequality holds in Lemma \ref{lemma:CompareRNR1} when $x > \lfloor x \rfloor$. Also, an argument similar to that in Lemmas \ref{lemmaWI} and \ref{lemmaWx} can be used to show that for $1 \leq i \leq j \leq \lfloor x \rfloor+1$, $\hat{z}_{i,j}^{(\lfloor x \rfloor+1,x)}$ is strictly decreasing in $j$ for $1 \leq j \leq \lceil x \rceil +1$, and $x$ when $x > 1$.

\subsection{The Nash equilibrium and its comparison with the $N$-case} \label{sec: Rcase:NashE}
As in the $N$-case, to work out the Nash equilibrium in the $R$-case, we need to draw the best response plot and investigate the intersection point of $\mathcal{BR}(x)$ and $x$. When $\hat{z}^{(m+1, m)}_{m+1,m+1} \leq 0 \leq \hat{z}^{(m,m)}_{m,m}$, the tagged customer's best response when others adopt $m$ is also $m$, which is the case in Figure \ref{fig:BR}(a). When $\hat{z}^{(m+1, x)}_{m+1,m+1} = 0$ with $x \in (m,m+1)$, the tagged customer is indifferent between $m$ and $m+1$ when others use threshold $x$, which is the case in Figure \ref{fig:BR}(b). 

\begin{figure}[h]
	\centering{\includegraphics[width=0.36\linewidth]{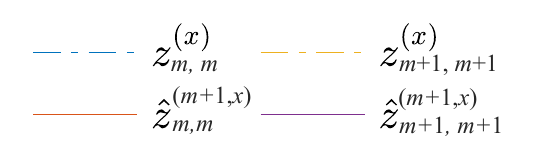}}\\
	\subcaptionbox{$R_0 \in [\alpha_m, \gamma_m]$}
	{\includegraphics[width=0.32\linewidth]{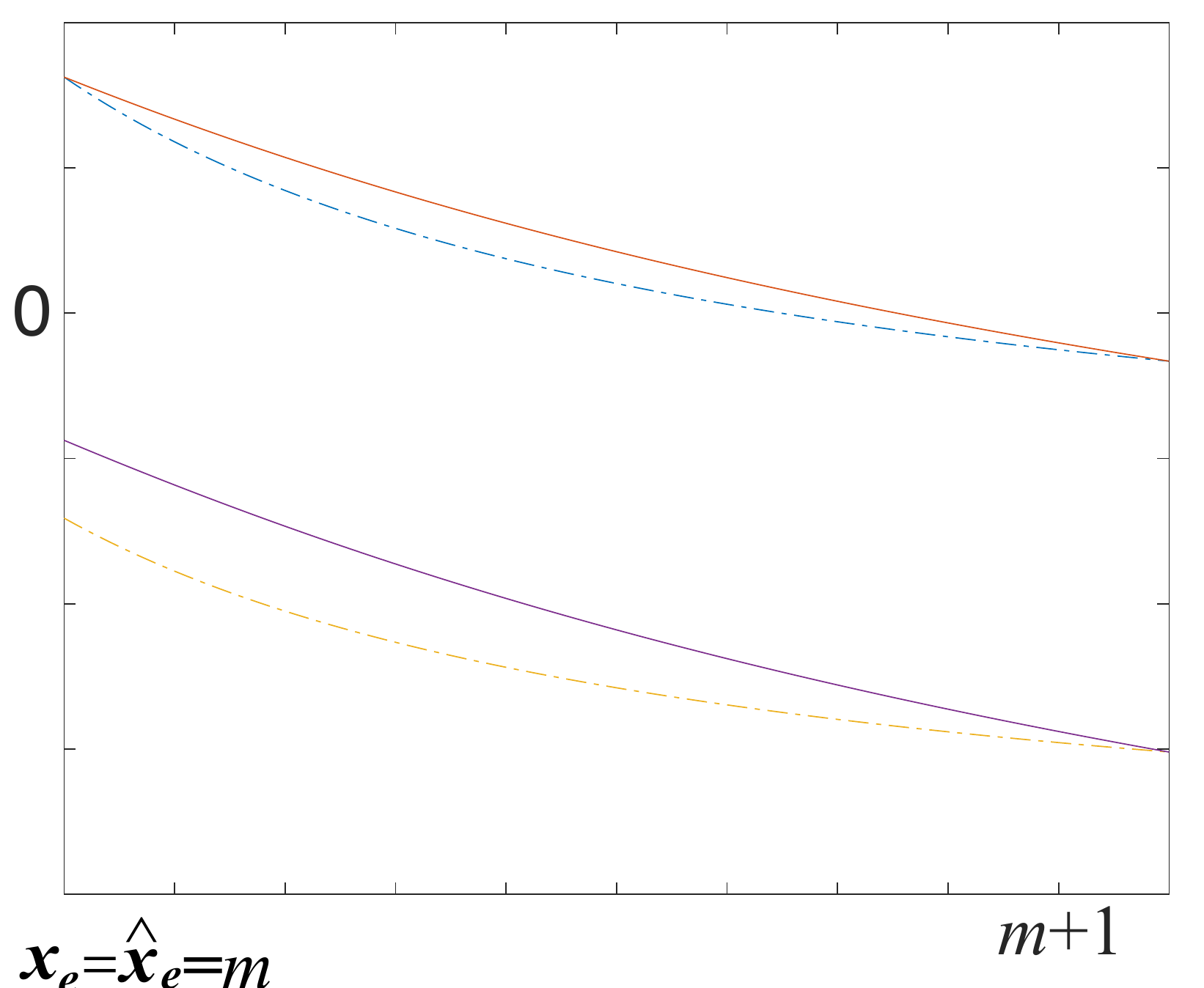}}  
	\hspace{\fill}
	\subcaptionbox{$R_0 \in (\gamma_m, \beta_m]$}
	{\includegraphics[width=0.32\linewidth]{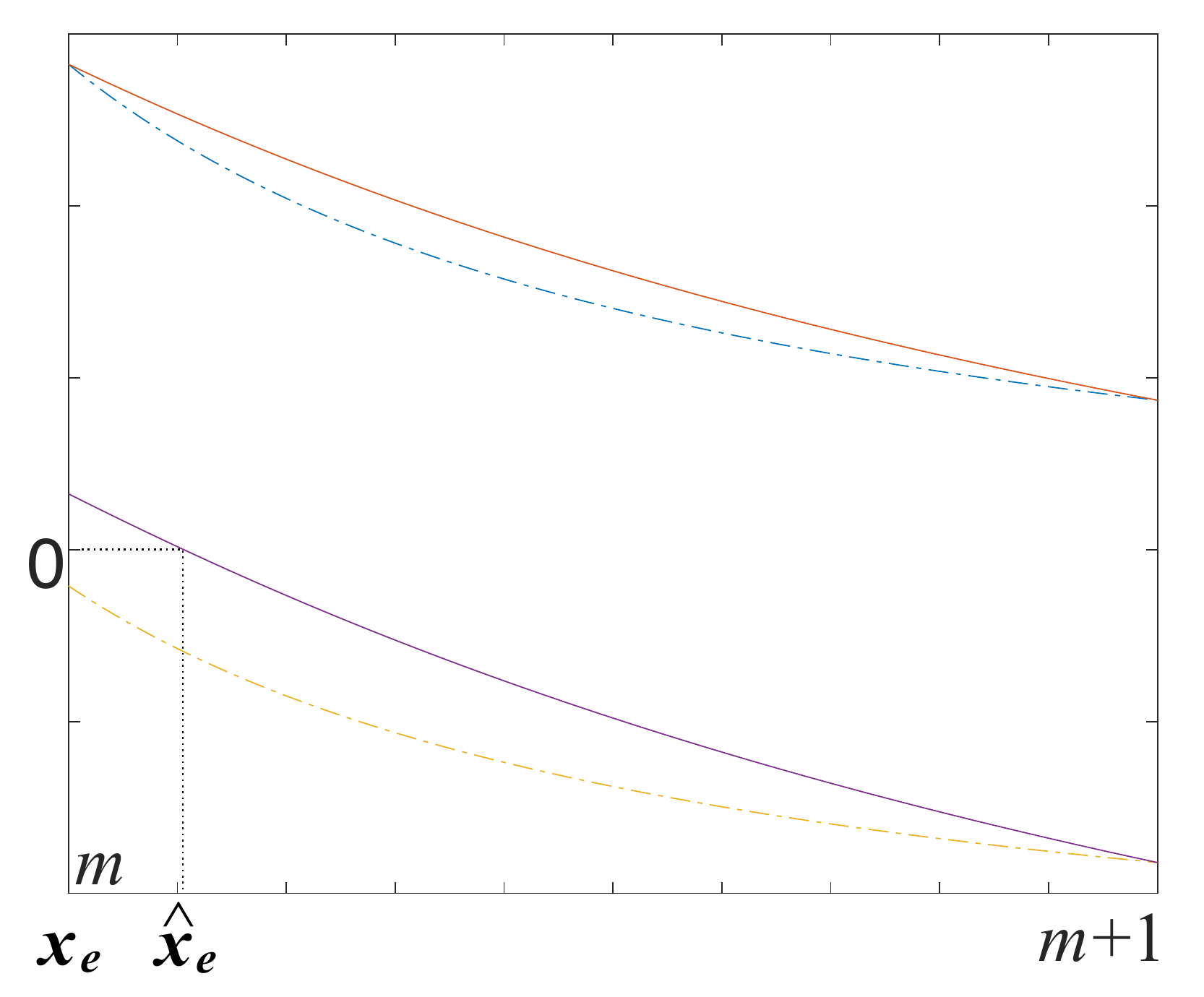}} 
	\hspace{\fill}
	\subcaptionbox{$R_0 \in (\beta_m, \alpha_{m+1})$}
	{\includegraphics[width=0.32\linewidth]{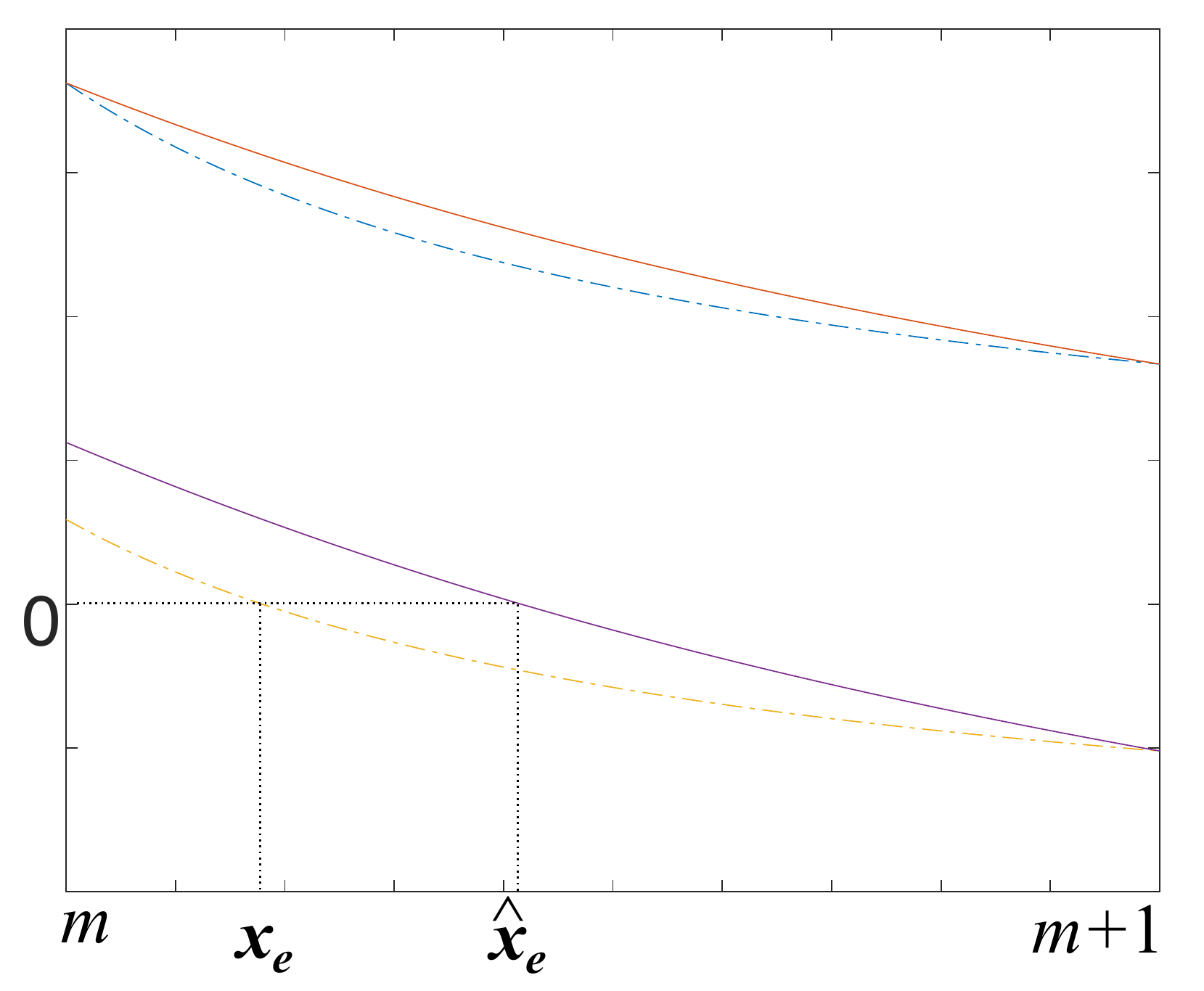}} 
	\caption{An illustration of Nash equilibrium threshold comparison.} \label{NER}
\end{figure}

Before we work out the Nash equilibrium strategy in the $R$-case and compare it with the $N$-case, we first define $Ne(R_0, \lambda, \mu, q)$ and $\hat{Ne}(R_0, \lambda, \mu, q)$ as the Nash equilibrium under the parameter set $R_0, \lambda, \mu, q$ in the $N$-case and the $R$-case, respectively, and use $x_e$ and $\hat{x}_e$ for short if they are from the same $R_0, \lambda, \mu, q$. Similar to our use of $\alpha_m$ and $\beta_m$ in the $N$-case, we let $\gamma_m := \hat{w}_{m+1,m+1}^{(m+1,m)}$ to help explain the Nash equilibrium in the $R$-case which is described in the following. 
\begin{theorem} \label{theorem:NER}
	The Nash equilibrium threshold value when reneging is allowed is greater than or equal to that when reneging is not allowed.
\end{theorem}
Proof.  There are three scenarios. 
\begin{itemize}
	\item When $R_0 \in [\alpha_m, \gamma_m]$, then
	\begin{align}
	& \hat{z}_{m+1,m+1}^{(m+1,m)}= R_0 - \hat{w}_{m+1,m+1}^{(m+1,m)} = R_0-\gamma_m \leq 0 \\
	& z_{m,m}^{(m)} = R_0 - z_{m,m}^{(m)} = R_0-\alpha_m \geq 0 \,.
	\end{align} 
	Hence $z_{m+1,m+1}^{(m)} < \hat{z}_{m+1,m+1}^{(m+1,m)} \leq 0 \leq  \hat{z}_{m,m}^{(m,m)} = z_{m,m}^{(m)}$ with the first inequality and the equality following from Lemma \ref{lemma:CompareRNR1}. The tagged customer's best response is $m$ if others' strategy is $m$ in both the $N$-case and the $R$-case, hence $\hat{x}_e = x_e = m$. This case is depicted in Figure \ref{NER}(a).
	\item When $R_0 \in (\gamma_m, \beta_m]$, then
	\begin{align}
	& \hat{z}_{m+1,m+1}^{(m+1,m)}= R_0 - \hat{w}_{m+1,m+1}^{(m+1,m)} = R_0-\gamma_m > 0 \\
	& z_{m+1,m+1}^{(m)} = R_0 - z_{m+1,m+1}^{(m)} = R_0-\beta_m \leq 0 \,.
	\end{align}  
	Hence $\hat{z}_{m+1,m+1}^{(m+1,m+1)}  < z_{m+1,m+1}^{(m)} \leq 0 < \hat{z}_{m+1,m+1}^{(m+1,m)} <  \hat{z}_{m,m}^{(m+1,m)} = z_{m,m}^{(m)}$ with the first inequality, the last inequality and the equality following from Lemma \ref{lemma:CompareRNR1}. The tagged customer's best response is $m$ if others' strategy is $m$ in the $N$-case. In the $R$-case, $\hat{x}_e = \lbrace x: \, \hat{z}_{m+1,m+1}^{(m+1, x)} = 0 \rbrace$ since the tagged customer is indifferent between joining at position $m+1$ and balking if others adopt threshold $\hat{x}_e$. In this case, $\hat{x}_e > x_e = m$, and it is depicted in Figure \ref{NER}(b).
	\item When $R_0 \in (\beta_m, \alpha_{m+1})$, then
	\begin{equation}
	\hat{z}_{m+1,m+1}^{(m+1,m+1)} = z_{m+1,m+1}^{(m+1)} < 0 < z_{m+1,m+1}^{(m)} < \hat{z}_{m+1,m+1}^{(m+1,m)} \,,
	\end{equation}
	with the first equality and the last inequality following from Lemma \ref{lemma:CompareRNR1}. The Nash equilibrium strategies are
	\[
	x_e = \lbrace x : \, z_{m+1,m+1}^{(x)} = 0 \rbrace \qquad \hat{x}_e = \lbrace x : \hat{z}_{m+1,m+1}^{(m+1, x)} = 0 \rbrace \,,
	\]
	which are mixed in both cases. It follows from Equation \eqref{eq:compare1} that $\hat{x}_e > x_e$. Note that $x_e$ here is the same as calculated in Theorem \ref{Theorem:NE}.  This case is depicted in Figure \ref{NER}(c). \hfill $\square$
\end{itemize}

\section{Two Paradoxes} \label{sec:Paradox}
In the $N$-case, every customer remains in the system until she successfully completes her service and receives reward $R_0$. Increasing $R_0$ can increase customers' incentive to join but also make the system more crowded. In this situation, does everyone become better off when the reward $R_0$ increases? To answer this question, we observe that there are parameter settings where the equilibrium expected payoff can decrease with $R_0$.  This paradoxical behaviour is discussed in the following.
\begin{paradox} \label{paradox1}
	In the $N$-case, let $x_k = NE(r_k, \lambda, \mu, q), \, k= 1,2.$ Then for $m = 1,2$, $z^{(x_1)}_{m,m} > z^{(x_2)}_{m,m}$ if $\beta_m < r_1 < r_2 < \alpha_{m+1}$, where $m$ is the integer part of the Nash equilibrium. In other words, if $R_0 \in (\beta_m, \alpha_{m+1})$, increasing $R_0$ will make everyone joining at position $m$ worse off.
\end{paradox}
Proof. As in Definition \ref{D1}, $p$ is the fractional part of $x$. When $x \in (1,2)$, 
\begin{equation}
w_{2,2}^{(x)} - w_{1,1}^{(x)} = \frac{\mu+\lambda p}{\mu (-\mu q^2 + 2\mu q + \lambda p)} = \frac{1}{\mu \, \left(1-\frac{\mu\,(1-q)^2}{\mu+\lambda \, p}\right)} \,,
\end{equation}
which is decreasing in $p$. When $x \in (2,3)$, 
\begin{equation}
w_{3,3}^{(x)} - w_{2,2}^{(x)} = \frac{1}{\mu (1-\mu^2 (1-q)^2\, f(p))} \,,
\end{equation}
where $\displaystyle f(p)=\frac{\lambda +2 \lambda  p q+\mu  q^3-3 \mu  q^2-\lambda  q+3 \mu  q}{(\mu +\lambda  p) \left(\lambda  \mu +\lambda ^2 p+\lambda  \mu  p q+\mu ^2 q^3-3 \mu ^2 q^2+3 \mu ^2 q\right)}$, which is decreasing in $p$. See Appendix \ref{appendix:fp} for the derivative of function $f(p)$. Hence $w_{3,3}^{(x)} - w_{2,2}^{(x)} $ is decreasing in $p$.

From Theorem \ref{Theorem:NE}, if $\beta_m < r_1 < r_2 < \alpha_{m+1}$, then $x_k$, which is the Nash equilibrium threshold when $R= r_k$, satisfies
\begin{equation}
w_{m+1,m+1}^{(x_k)} = r_k \qquad x_k \in (m,m+1) \qquad k = 1,2 \,.
\end{equation} 
Thus, for $m =1,2$,
\begin{equation}
\Scale[0.9]{
z_{m,m}^{(x_1)} = r_1 - w_{m,m}^{(x_1)} =w_{m+1,m+1}^{(x_1)}   - w_{m,m}^{(x_1)} >  w_{m+1,m+1}^{(x_2)} - w_{m,m}^{(x_2)} = r_2 - w_{m,m}^{(x_2)} = z_{m,m}^{(x_2)} \,. }
\end{equation} 
\hfill $\square$

In Paradox \ref{paradox1}, we have proved that if $R_0 \in (\beta_m, \alpha_{m+1})$ for $m = 1,2$, increasing $R_0$ makes everyone joining at position $m$ worse off. We conjecture that this phenomenon holds for any $m$. Our numerical experience indicates that this is the case. However, the proof has eluded us. 

We have proved in the previous section that customers have a higher incentive to join the system if they are allowed to renege later. However, with more customers joining, the system can be more crowded. So we are interested in the question: when customers are given the right to leave, do they become better off? 

To answer this question, we first need to work out the equilibrium expected payoff in the $R$-case. 
In contrast to the $N$-case, customers may renege before they successfully complete the service in the $R$-case, thus their expected payoff cannot be calculated as the difference between $R_0$ and their expected sojourn cost. 
By similar reasoning to Equation \eqref{eq:wij}, it follows that 
\begin{align}
&\Scale[0.97]{\hat{z}_{i,j}^{(x,x)} =  -\frac{1}{\lambda+\mu}}\\
&\Scale[0.97]{+ \frac{\lambda}{\lambda+\mu} \, \left( \hat{z}_{i,j+1}^{(x,x)} \, \mathbbm{1}_{\lbrace j < \lfloor x \rfloor \rbrace} + \left(p \, \hat{z}_{i,j+1}^{(x,x)} + (1-p) \, \hat{z}_{i,j}^{(x,x)}\right) \, \mathbbm{1}_{\lbrace j = \lfloor x \rfloor \rbrace} +\hat{z}_{i,j}^{(x,x)} \, \mathbbm{1}_{\lbrace j = \lfloor x \rfloor+1 \rbrace} \right) \notag } \\
& \,\Scale[0.97]{+ \frac{\mu}{\lambda+\mu} \, \left[\left( q \, R_0 + (1-q) \, (1-(1-p) \mathbbm{1}_{\lbrace j=\lfloor x \rfloor+1 \rbrace}) \, \hat{z}_{j,j}^{(x,x)} \right) \, \mathbbm{1}_{\lbrace  i = 1 \rbrace} \right. }\\
&\Scale[0.97]{\left. \, + \left( (q + (1-q) \, (1-p) \, \mathbbm{1}_{\lbrace j=\lfloor x \rfloor+1 \rbrace}) \, \hat{z}_{i-1,j-1}^{(x,x)} + (1-q) \, (1-(1-p) \mathbbm{1}_{\lbrace j=\lfloor x \rfloor+1 \rbrace}) \, \hat{z}_{i-1,j}^{(x,x)}  \right) \, \mathbbm{1}_{\lbrace i > 1 \rbrace}\right] \notag .}
\end{align}
Thus, we can obtain $\hat{z}_{i,j}^{(x,x)} $ via Poisson's equation
\begin{equation}
\left(I-\hat{P}^{(x,x)} \right) \, 
\bm{\hat{z}}^{(x,x)} \, 
=  \, \bm{{g}} \,,
\end{equation}
where $\hat{P}^{(x,x)}$ is defined in Appendix \ref{appendix: 3.1}, and
\begin{align}
\begin{split}
&\bm{\hat{z}}^{(x,x)} = (z_{1,1}^{(x,x)}, z_{1,2}^{(x,x)}, z_{2,2}^{(x,x)}, z_{1,3}^{(x,x)}, z_{2,3}^{(x,x)}, z_{3,3}^{(x,x)},\cdots, z_{\lceil x \rceil,\lceil x \rceil+1}^{(x,x)}, z _{\lceil x \rceil+1,\lceil x \rceil+1}^{(x,x)})^T \,.
\end{split}
\end{align}
We are interested in the expected payoff when every customer uses $\hat{x}_e$. In Lemma \ref{lemma:CompareRNR1}, we have proved that if $\hat{x}_e = \lfloor \hat{x}_e \rfloor$, $\bm{\hat{z}}^{(\hat{x}_e,\hat{x}_e)}=\bm{{z}}^{(x_e)}$. In the following, we prove that if $\hat{z}^{(m+1,\hat{x}_e)}_{m+1,m+1} = 0$, $\bm{\hat{z}}^{(\hat{x}_e,\hat{x}_e)} = \bm{\hat{z}}^{(m+1,\hat{x}_e)}$.
\begin{lemma} \label{lemma:mm1same}
	If $\hat{z}^{(m+1,\hat{x}_e)}_{m+1,m+1} = 0$, then $\hat{z}^{(\hat{x}_e,\hat{x}_e)}_{i,j} = \hat{z}^{(m+1,x_e)}_{i,j}$ for any $1 \leq i \leq j \leq m+1$.
\end{lemma}
Proof.  First, when $\hat{z}^{(m+1,\hat{x}_e)}_{m+1,m+1} = 0$, the tagged customer is indifferent between joining or not joining at position $m+1$. In other words, joining with any probability at position $m+1$ will result in a zero expected payoff for her. Hence, $\hat{z}^{(x,\hat{x}_e)}_{m+1,m+1} = 0$ for any $x  \in [m, m+1]$ including $\hat{x}_e$.

For a general state $(i,j)$, consider two queues with others using threshold $\hat{x}_e$ and the tagged customer in state $(i,j)$: she uses threshold $m+1$ in queue 1 and $\hat{x}_e$ in queue 2.  By coupling the customer arrival processes, their joining decisions, the service processes and the service success probability for every customer including the tagged one, we can see the next customer will arrive, join or not join both queues at the same time, the customer in service in both queues will complete the service and rejoin or not rejoin the queue at the same time, until the first time the tagged customer needs to rejoin the queue and the queue size including her is $m+1$ when she rejoins it. When this is the case, $\hat{z}^{(x, \hat{x}_e)}_{m+1,m+1} = 0$ for any $x \in [m,m+1]$, due to $\hat{z}^{(m+1,\hat{x}_e)}_{m+1,m+1} = 0$. Hence, the tagged customer in both queues either has exactly the same sample path, or reaches the state $(m+1, m+1)$ where her remaining expected payoff is $0$ regardless of her decision. So the tagged customer in both queues receives the same expected payoff. \hfill $\square$

When $R_0 \in [\alpha_m, \gamma_m]$, then $x_e = \hat{x}_e=m$, thus, $z^{(m)}_{i,i} = \hat{z}^{(m,m)}_{i,i}$, and the stationary distributions in the $N$-case and the $R$-case are the same. Thus there is no difference between the two cases. When $R_0 \in (\gamma_m, \alpha_{m+1}]$, we observe that both the equilibrium expected payoff and the total expected payoff 
decrease when reneging is allowed.  Specifically, when $R_0 \in (\gamma_m, \beta_m)$, we prove this paradoxical behaviour in the following.
\begin{paradox} \label{paradox2}
	If $R_0 \in  (\gamma_m, \beta_m]$, then
	\[
	\hat{z}_{i,i}^{(\hat{x}_e,\hat{x}_e)} < z_{i,i}^{({x}_e)} \quad i = 1, \cdots, m \,.
	\]Furthermore, the total expected payoff under equilibrium satisfy
	\[
	\sum_{i = 1}^{m+1} \hat{\pi}^{(\hat{x}_e)}_{i-1} \, \hat{z}_{i,i}^{(\hat{x}_e,\hat{x}_e)} < \sum_{i = 1}^{m}\pi^{(x_e)}_{i-1} z_{i,i}^{({x}_e)}
	\]
	where $\pi^{(x)}_k$ and $\hat{\pi}^{(x)}_k$ denote the stationary distribution of the number of customers in the system in the $N$-case and the $R$-case, respectively. 
\end{paradox}
Proof.
If $R_0 \in  (\gamma_m, \beta_m]$, then $\hat{x_e} > x_e = m$. Since the Nash equilibrium threshold in the $R$-case is mixed, it follows from Lemma \ref{lemma:mm1same} that $z_{m+1,m+1}^{(\hat{x}_e, \hat{x}_e)} = z_{m+1,m+1}^{(m+1, \hat{x}_e)} = 0$, and $z_{i,i}^{(\hat{x}_e, \hat{x}_e)} = z_{i,i}^{(m+1, \hat{x}_e)} $. In this way, 
\begin{equation} \label{eq:com1}
z^{(x_e)}_{i,i}  = z^{(m)}_{i,i}= \hat{z}^{(m+1, m)}_{i,i} > \hat{z}^{(m+1, \hat{x}_e)}_{i,i} = \hat{z}^{(\hat{x}_e, \hat{x}_e)}_{i,i} \qquad 1 \leq i \leq m\,,
\end{equation}
with the second equality following from Lemma \ref{lemma:CompareRNR1}, and the inequality following  from the fact that $\hat{z}_{m+1,m+1}^{(m+1,x)}$ is decreasing in $x$.

\begin{figure}[H]
	\centering
	\begin{tikzpicture}[line cap=round,line join=round,>=triangle 45,x=1.0cm,y=1.0cm, scale=0.6 , transform shape] 
	
	\clip(-2.84,-1.3) rectangle (9.04,3.58);
	\node (zero) [draw, circle, minimum size=1.4cm] at (-2,1cm) {$0$};
	\node (one) [draw, circle, minimum size=1.4cm] at (0,1cm) {$1$};
	\node (two) [draw, circle, minimum size=1.4cm] at (2,1cm) {$2$};
	\node (enminusone) [draw, circle, minimum size=1cm] at (4,1cm) {$\lfloor x \rfloor-1$};
	\node (en) [draw, circle, minimum size=1.4cm] at (6,1cm) {$\lfloor x \rfloor$};
	\node (enplusone) [draw, circle, minimum size=1cm] at (8,1cm) {$\lfloor x \rfloor+1$};
	
	\draw [->] (zero) .. controls +(0.5,1.5) and +(-0.5,1.5) .. node [midway, above] {$\lambda$} (one);
	\draw [->] (one) .. controls +(0.5,1.5) and +(-0.5,1.5) .. node [midway, above] {$\lambda$} (two);
	\draw [->, very thick, dashed] (two) .. controls +(0.5,1.5) and +(-0.5,1.5) .. (enminusone);
	\draw [->] (enminusone) .. controls +(0.5,1.5) and +(-0.5,1.5) .. node [midway, above] {$\lambda$} (en);
	\draw [->] (en) .. controls +(0.5,1.5) and +(-0.5,1.5) .. node [midway, above] {$\lambda p$} (enplusone);
	
	\draw [->] (enplusone) .. controls +(-0.5,-1.5) and +(0.5,-1.5) .. node [midway, below] {$\mu q$} (en);
	\draw [->] (en) .. controls +(-0.5,-1.5) and +(0.5,-1.5) .. node [midway, below] {$\mu q$} (enminusone);
	\draw [->, very thick, dashed] (enminusone) .. controls +(-0.5,-1.5) and +(0.5,-1.5) .. (two);
	\draw [->] (two) .. controls +(-0.5,-1.5) and +(0.5,-1.5) .. node [midway, below] {$\mu q$} (one);
	\draw [->] (one) .. controls +(-0.5,-1.5) and +(0.5,-1.5) .. node [midway, below] {$\mu q$} (zero);
	\end{tikzpicture}
	\caption{Transition rate diagram when the threshold is $x$ when reneging is not allowed.} \centering \label{fig:TR1}
	\begin{tikzpicture}[line cap=round,line join=round,>=triangle 45,x=1.0cm,y=1.0cm,scale=0.6 , transform shape]
	
	\clip(-2.84,-1.3) rectangle (9.04,3.58);
	\node (zero) [draw, circle, minimum size=1.4cm] at (-2,1cm) {$0$};
	\node (one) [draw, circle, minimum size=1.4cm] at (0,1cm) {$1$};
	\node (two) [draw, circle, minimum size=1.4cm] at (2,1cm) {$2$};
	\node (enminusone) [draw, circle, minimum size=1cm] at (4,1cm) {$\lfloor x \rfloor-1$};
	\node (en) [draw, circle, minimum size=1.4cm] at (6,1cm) {$\lfloor x \rfloor$};
	\node (enplusone) [draw, circle, minimum size=1cm] at (8,1cm) {$\lfloor x \rfloor+1$};
	
	\draw [->] (zero) .. controls +(0.5,1.5) and +(-0.5,1.5) .. node [midway, above] {$\lambda$} (one);
	\draw [->] (one) .. controls +(0.5,1.5) and +(-0.5,1.5) .. node [midway, above] {$\lambda$} (two);
	\draw [->, very thick, dashed] (two) .. controls +(0.5,1.5) and +(-0.5,1.5) .. (enminusone);
	\draw [->] (enminusone) .. controls +(0.5,1.5) and +(-0.5,1.5) .. node [midway, above] {$\lambda$} (en);
	\draw [->] (en) .. controls +(0.5,1.5) and +(-0.5,1.5) .. node [midway, above] {$\lambda p$} (enplusone);
	
	\draw [->] (enplusone) .. controls +(-0.5,-1.5) and +(0.5,-1.5) .. node [midway, below] {\small $\mu q + \mu (1-q)(1-p)$} (en);
	\draw [->] (en) .. controls +(-0.5,-1.5) and +(0.5,-1.5) .. node [midway, below] {$\mu q$} (enminusone);
	\draw [->, very thick, dashed] (enminusone) .. controls +(-0.5,-1.5) and +(0.5,-1.5) .. (two);
	\draw [->] (two) .. controls +(-0.5,-1.5) and +(0.5,-1.5) .. node [midway, below] {$\mu q$} (one);
	\draw [->] (one) .. controls +(-0.5,-1.5) and +(0.5,-1.5) .. node [midway, below] {$\mu q$} (zero);
	\end{tikzpicture}
	\caption{Transition rate diagram when the threshold is $x$ when reneging is allowed.} \centering \label{fig:TR2}
\end{figure}

Next, we calculate the stationary distribution of the number of customers in the system in the $N$-case and the $R$-case. Figures \ref{fig:TR1} and \ref{fig:TR2} depict the transition rate diagram for both cases, given that each customer uses threshold $x$. Let $\displaystyle \rho = \frac{\lambda}{\mu q}$, it follows from the detailed balance equations that for $k = 0, \cdots, \lfloor x \rfloor$,
\begin{align}
&\pi_k^{(x)}  = \frac{\displaystyle\rho^{k}}{\displaystyle(x-\lfloor x \rfloor) \,\rho^{\lfloor x \rfloor+1}+\frac{\rho^{\lfloor x \rfloor+1}-1}{\rho-1}}  \\
&\hat{\pi}_k^{(x)}  = \frac{\displaystyle\rho^{k}}{\displaystyle\frac{\lambda (x - \lfloor x \rfloor)}{\mu q+\mu(1-q)(1-(x - \lfloor x \rfloor))} \rho^{\lfloor x \rfloor}+\frac{\rho^{\lfloor x \rfloor+1}-1}{\rho-1}} \notag \,,
\end{align}
and
\begin{align} 
	&\pi_{ \lfloor x \rfloor
	+1}^{(x)} = \frac{\displaystyle(x-\lfloor x \rfloor) \, \rho^{n+1}}{\displaystyle(x - \lfloor x \rfloor) \,\rho^{\lfloor x \rfloor+1}+\frac{\rho^{\lfloor x \rfloor+1}-1}{\rho-1}} \\
 	&\hat{\pi}_{\lfloor x \rfloor+1}^{(x)} = \frac{\displaystyle \frac{\lambda \, (x - \lfloor x \rfloor)}{\mu q+\mu(1-q)(1-(x - \lfloor x \rfloor))} \rho^{\lfloor x \rfloor}}{\displaystyle \frac{\lambda \, (x - \lfloor x \rfloor)}{\mu q+\mu(1-q)(1-(x - \lfloor x \rfloor))} \rho^{\lfloor x \rfloor}+\frac{\rho^{\lfloor x \rfloor+1}-1}{\rho-1}} \,. \notag
\end{align}
Since $\hat{x}_e > x_e = m$,
\begin{equation} \label{eq:SdCom}
\hat{\pi}^{(\hat{x}_e)}_k < \pi^{(x_e)}_k \quad k=0, \cdots, m \qquad \hat{\pi}^{(\hat{x}_e)}_{m+1} > 0 = \pi^{(x_e)}_{m+1} \,.
\end{equation} 
Hence
\begin{equation}
\sum_{i = 1}^{m+1} \hat{\pi}^{(\hat{x}_e)}_{i-1} \, \hat{z}_{i,i}^{(\hat{x}_e,\hat{x}_e)} = \sum_{i = 1}^{m} \hat{\pi}^{(\hat{x}_e)}_{i-1} \, \hat{z}_{i,i}^{(\hat{x}_e,\hat{x}_e)} < \sum_{i = 1}^{m}\pi^{(x_e)}_{i-1} z_{i,i}^{({x}_e)} \,,
\end{equation}
with the equality following from $z_{m+1,m+1}^{(m+1, \hat{x}_e)} = 0$, and the inequality following from Equations \eqref{eq:com1} and \eqref{eq:SdCom}.
\hfill $\square$

It can be seen that although the expected payoff $z_{m+1,m+1}^{(x_e)}$ is smaller than $\hat{z}_{m+1,m+1}^{(\hat{x}_e,\hat{x}_e)}$, customers do not really join at position $m+1$ in the $R$-case equilibrium, thus it is not included in the total expected payoff. 

\begin{figure}[H]
	\centering{\includegraphics[width=0.38\linewidth]{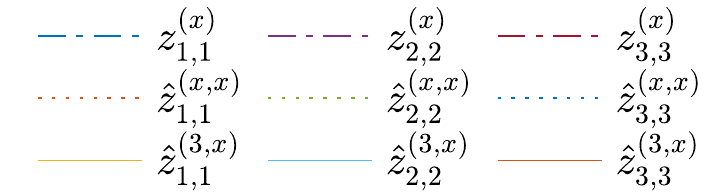}}\\
	\subcaptionbox{$R_0 \in (\gamma_m, \beta_m]$}%
	{\includegraphics[width=0.5\linewidth]{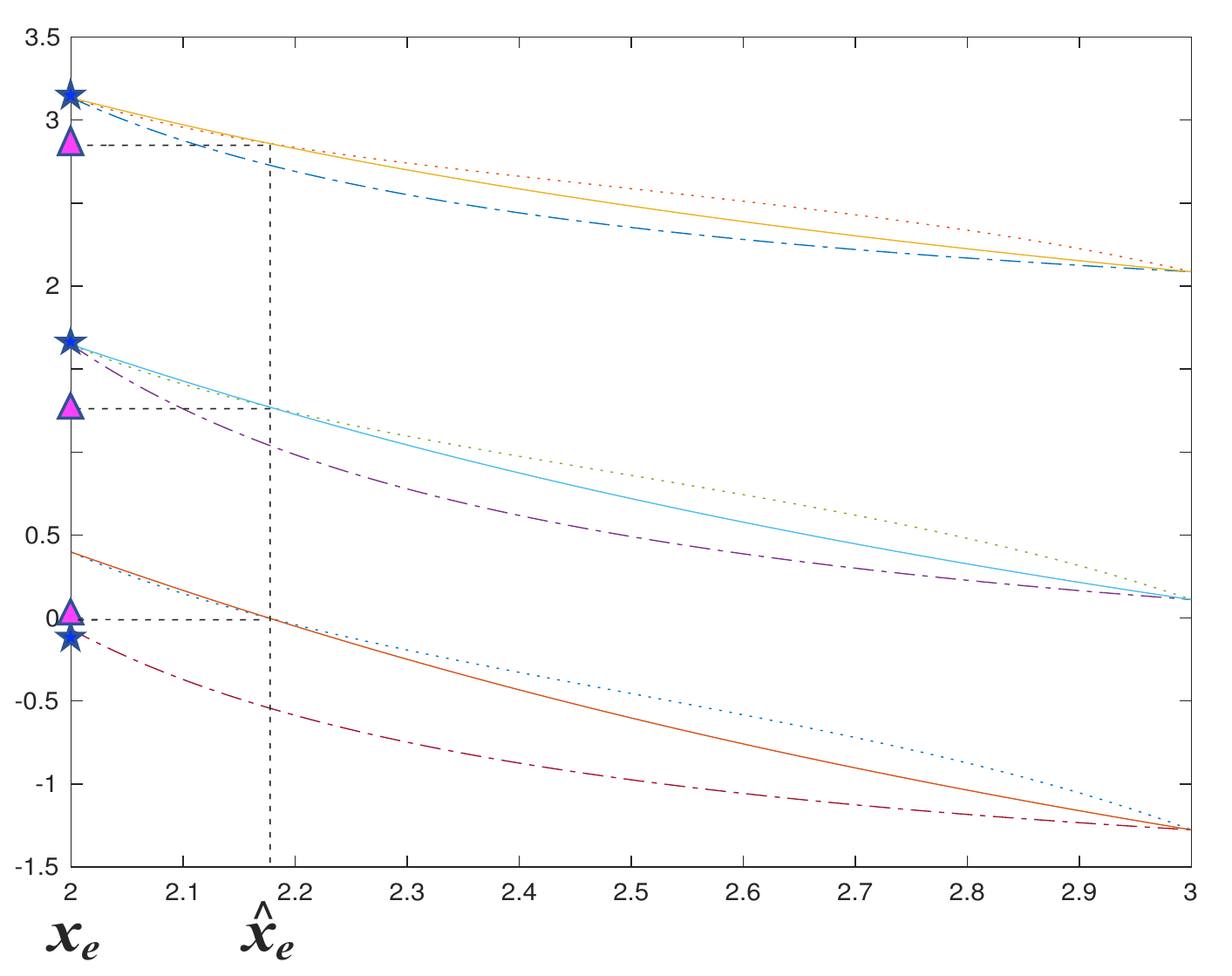}}
	\subcaptionbox{$R_0 \in (\beta_m, \alpha_{m+1})$}%
	{\includegraphics[trim={0 1.5mm  0  0 },clip, width=0.5\linewidth]{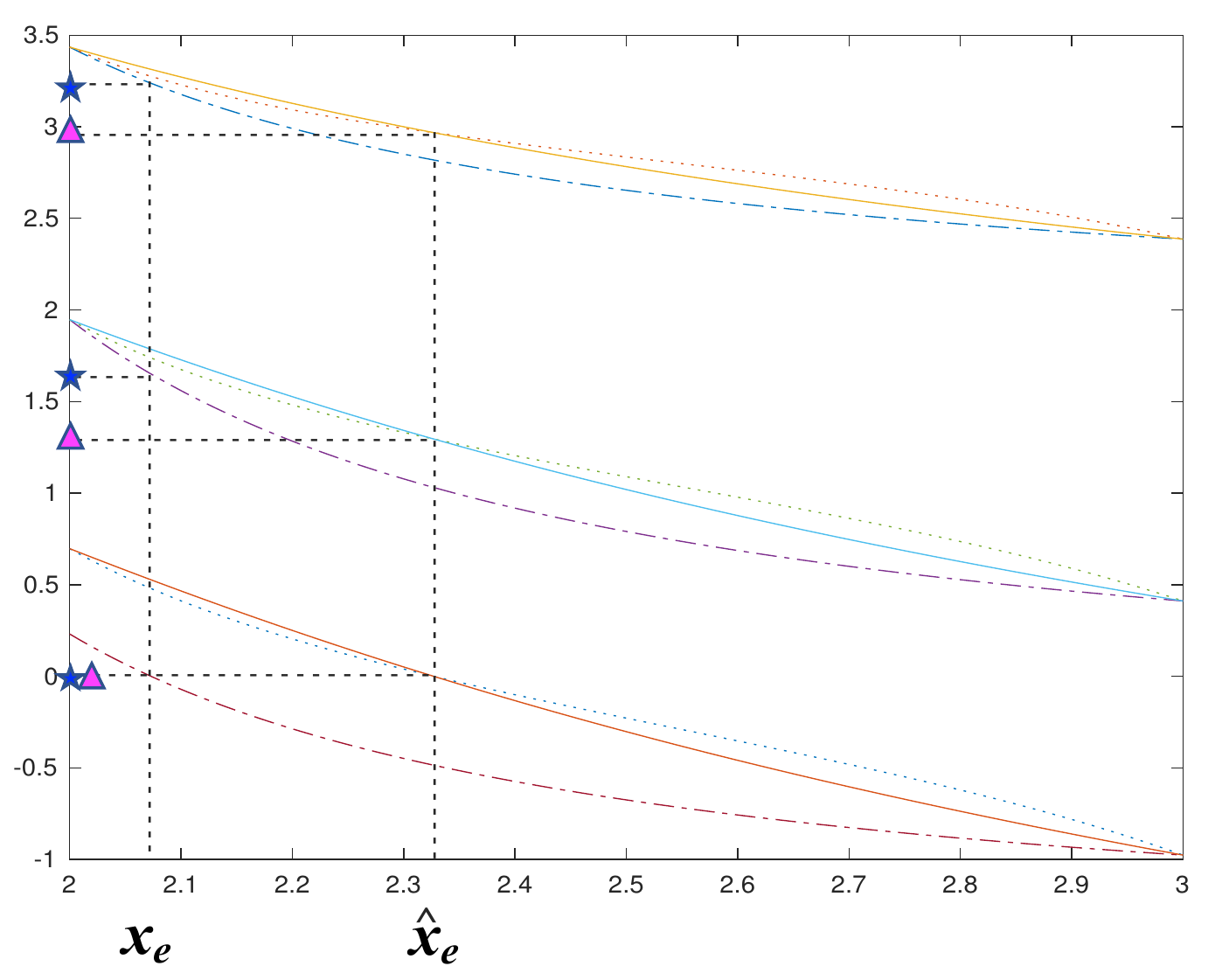}} 
	\caption{Allowing reneging can make everyone worse off.} \label{fig:pa2}
\end{figure}

We illustrate Equation (\ref{eq:com1}) in Figure \ref{fig:pa2}(a) via an example with $R_0 = 7.5, \lambda = 1, \mu = 0.8, q = 0.4$. The Nash equilibrium threshold is $2$ and $2.167$ in the $N$-case and the $R$-case, respectively. The blue stars represent $\bm{z}^{(x_e)}$ with $z^{(x_e)}_{1,1} > z^{(x_e)}_{2,2} > 0 > z^{(x_e)}_{3,3}$, and the triangles represent $\bm{\hat{z}}^{(\hat{x}_e, \hat{x}_e)}$ with $z^{(\hat{x}_e, \hat{x}_e)}_{1,1} > z^{(\hat{x}_e, \hat{x}_e)}_{2,2} > z^{(\hat{x}_e, \hat{x}_e)}_{3,3} = 0$. It can be observed that and $z_{i,i}^{(x_e)} > \hat{z}_{i,i}^{(\hat{x}_e,\hat{x}_e)}$ for $i =1,2$. 

In Paradox \ref{paradox2}, we proved that when $R_0 \in (\gamma_m, \beta_m]$, allowing reneging makes everyone worse off. Next, we use some numerical examples to show that allowing reneging can make everyone worse off when $R_0 \in (\beta_m, \alpha_{m+1})$.
Actually, the paradox is observed in every example. 

We first illustrate $z_{i,i}^{(x_e)}>\hat{z}_{i,i}^{(x_e,x_e)}$ 
for any $1 \leq i \leq \lfloor x_e \rfloor$ in Figure \ref{fig:pa2}(b) via an example with $R_0 = 7.8, \lambda = 1, \mu = 0.8, q = 0.4$. The Nash equilibrium is $2.073$ and $2.327$ in the $N$-case and the $R$-case, respectively. The blue stars represent $\bm{z}^{(x_e)}$ with $z^{(x_e)}_{1,1} > z^{(x_e)}_{2,2} > z^{(x_e)}_{3,3} = 0$, and the triangles represent $\bm{\hat{z}}^{(\hat{x}_e, \hat{x}_e)}$ with $z^{(\hat{x}_e, \hat{x}_e)}_{1,1} > z^{(\hat{x}_e, \hat{x}_e)}_{2,2} > z^{(\hat{x}_e, \hat{x}_e)}_{3,3} = 0$. It can be observed that $z_{i,i}^{(x_e)} > \hat{z}_{i,i}^{(\hat{x}_e,\hat{x}_e)}$ for $i =1,2$. 


\begin{table}
	\centering
	\begin{tabular}{|l | c | c  | c | } 
		\hline
		\quad $R_0,\lambda, \mu,q$ & $7.8, 1, 0.8, 0.4 $ & $4.4, 1, 0.8, 0.8$ & $13.5, 0.8, 1, 0.2 $ \, \\ [0.5ex] 
		\hline 
		$x_e \qquad \, \hat{x}_e$ & $2.073 \quad 2.327$ & $2.345 \quad 2.444$ & $2.529 \quad 2.872$ \\
		\hline
		$z_{1,1}^{(x_e)} \quad \hat{z}_{1,1}^{(\hat{x}_e, \hat{x}_e)}$ & $ 3.245 \quad 2.964$ & $2.599 \quad 2.591$ & $3.740 \quad 3.546$ \\
		\hline
		$z_{2,2}^{(x_e)} \quad \hat{z}_{2,2}^{(\hat{x}_e, \hat{x}_e)}$ & $1.661 \quad 1.292$ & $1.271 \quad 1.259$ & $1.514 \quad 1.283$ \\
		\hline
		$\pi_{0}^{(x_e)} \quad \hat{\pi}_{0}^{(\hat{x}_e)}$ & $0.063 \quad 0.053$ & $0.158 \quad 0.154$ & $ 0.018 \quad 0.017$ \\
		\hline
		$\pi_{1}^{(x_e)} \quad \hat{\pi}_{1}^{(\hat{x}_e)}$ & $0.195 \quad 0.165$ & $0.247 \quad 0.241$ & $0.073 \quad 0.069 $ \\
		\hline
\end{tabular}
\caption{An numerical example when $R_0 \in (\beta_m, \alpha_{m+1})$} \label{tab:pa2}
\end{table}	

Next, we list the Nash equilibrium, the stationary probabilities of having $0$ or $1$ customers in the queue, and the expected payoff of three examples with different values of $R_0, \lambda, \mu$, and $q$ given in Table \ref{tab:pa2}, to show that not only $z_{i,i}^{(x_e)}>\hat{z}_{i,i}^{(\hat{x}_e,\hat{x}_e)}$, but also $\pi_{i-1}^{(x_e)}>\hat{\pi}_{i-1}^{(\hat{x}_e)}$ for any $1 \leq i \leq \lfloor x \rfloor$. It follows from the transition rate diagram in Figure \ref{fig:TR1} and \ref{fig:TR2} that, for $i =1, \ldots m-1$, the transition rates going from state $i$ to state $i+1$, and vice versa, are identical in the $N$-case and the $R$-case, so the only difference in the stationary distribution for the two cases is the normalisation constant. The greater constant in the $R$-case makes the first $m$ states have less probability mass than the $N$-case, and it is only the final one that compensates. The first example has the same parameters as that in Figure \ref{fig:pa2}. When the Nash equilibrium is fractional, $z_{m+1,m+1}^{(x_e)} = \hat{z}_{m+1,m+1}^{(\hat{x}_e, \hat{x}_e)} = 0$ where $m$ is the integer part of the Nash equilibrium, so we omit this in the table. Also, to compare $\sum_{i = 1}^{m+1} {\pi}^{({x}_e)}_{i-1} \, {z}_{i,i}^{(x_e)}$ and $\sum_{i = 1}^{m+1} \hat{\pi}^{(\hat{x}_e)}_{i-1} \, \hat{z}_{i,i}^{(\hat{x}_e,\hat{x}_e)}$, we only need to calculate ${\pi}^{({x}_e)}_i$ and $\hat{\pi}^{(\hat{x}_e)}_i$ for $i = 0, \cdots, m-1$, so we omit ${\pi}^{({x}_e)}_k$ and $\hat{\pi}^{(\hat{x}_e)}_k$ for $k = m,m+1$. In Table \ref{tab:pa2}, the Nash equilibrium thresholds of the three examples are all fractional and have the same integer part, that is, $2$. We observe that $z_{i,i}^{(x_e)} > \hat{z}_{i,i}^{(\hat{x}_e,\hat{x}_e)}$ and $\pi_{i-1}^{(x_e)}>\hat{\pi}_{i-1}^{(\hat{x}_e)}$ for $i =1,2$. 

\section{Social Welfare} \label{sec:SW}

In the previous section, we showed that allowing reneging can make every customer worse off. If the goal is to maximise the social welfare which is defined as the
total expected net benefit of all customers, how does the reneging affect the social welfare? In this section, we calculate and compare the optimal threshold  from the social point of view in the $N$-case and the $R$-case.
\subsection{Social welfare in the $N$-case}
When the customers all adopt threshold $x$, the state transition rate diagram in the non-reneging case is shown in Figure \ref{fig:TR1}, the social welfare
\begin{align}
S^{N} (x) :&= \lambda \left(\sum_{k = 1}^{\lfloor x \rfloor} \pi_{k-1}^{(x)} z_{k,k}^{(x)}+  (x-\lfloor x \rfloor) \pi_{\lfloor x \rfloor}^{(x)} z_{\lfloor x \rfloor+1,\lfloor x \rfloor+1}^{(x)}\right)\\
& = \lambda R_0 \left(\sum_{k = 0}^{\lfloor x \rfloor-1} \pi_k^{(x)}+  (x-\lfloor x \rfloor) \pi_{\lfloor x \rfloor}^{(x)} \right)-\sum_{k = 0}^{\lceil x \rceil} \, k \, \pi_{k}^{(x)} \,.
\end{align}
where $\rho := \displaystyle \frac{\lambda}{\mu q}$, and the second equality follows from Little's law. The explicit expression is in Appendix \ref{appendix:sw}.

\begin{proposition} \label{pro1}
	Social welfare $S^{N}(x)$ is unimodal.
\end{proposition}
Proof. 
We first take the derivative of $S^{N}(x)$,
\begin{equation}\label{Sod1}
\frac{dS^{N}(x)}{dx}
= 
\Scale[0.91]{\begin{cases}
\displaystyle \frac{\rho^{\lfloor x \rfloor} \left(R_0 \lambda (\rho-1)^2 - \rho \,\left(1 - 2 \rho + \lfloor x \rfloor (1-\rho) + \rho^{\lfloor x \rfloor+2}\right)\right)}{(1 + \rho^{\lfloor x \rfloor+1} ((x-\lfloor x \rfloor)(1-\rho) - 1))^2}   \qquad & \text{when }x > \lfloor x \rfloor \\
\text{undefined} & \text{when }x = \lfloor x \rfloor \,.
\end{cases}}
\end{equation}
To see that $S^{N}(x)$ is unimodal, let
\begin{equation}
f(k) = \left(1 - 2 \rho + k (1-\rho) + \rho^{2 + k} \right) \quad k = 0,1,2,\cdots \,,
\end{equation}
and observe that the numerator in the first Equation of (6.4) can be written as 
\[
\rho^{\lfloor x \rfloor}\left(R_0 \lambda(\rho-1)^2 -\rho f(\lfloor x \rfloor) \right)\,.
\]
When $\rho \neq 1$
\begin{equation}
f(k+1) - f(k) = \rho^{2 + k} (\rho - 1) + (1-\rho) = (1-\rho)\left(1-\rho^{2 + k} \right) > 0 \,.
\end{equation}
We assume that $R_0 > \frac{\ds 1}{\ds \mu q}$ to avoid the trivial case where the reward is smaller than the expected service time even if a customer does not have to wait. Hence
\begin{equation}
f(0) = (-1+\rho)^2 < \frac{R_0 \lambda }{\rho}(1-\rho)^2
\end{equation}
and
\begin{equation}
\lim\limits_{n \rightarrow \infty} f(k) = \infty > \frac{R_0 \lambda }{\rho}(1-\rho)^2,
\end{equation}
and so,
\begin{equation}
\rho \, f(0) < \cdots < 	R_0 \lambda (1-\rho)^2< \cdots < \rho \,\lim\limits_{n \rightarrow \infty} f(k) \,.
\end{equation}
Thus, there exists an integer $n_N^{\star}$ such that $\dfrac{dS^{N}}{dx}$ is increasing when $ x  \leq n_N^{\star}$; is decreasing when $ x  > n_N^{\star}$. That is, $n_N^{\star}$ is the socially optimal threshold.
$ \hfill \square $

It can be observed that $n_N^{\star} = \lfloor \nu \rfloor$, where $\nu$ satisfies 
\begin{equation}
R_0 \mu q  -\nu = \frac{\rho}{(1-\rho)^2} \left( \nu (1-\rho)- 1+ \rho^{\nu} \right) \,.
\end{equation}
This coincides with Naor's result for the non-feedback $M/M/1$ queue in \citet[section 4]{N69}. In other words, from the perspective of society, the feedback parameter $q$ affects the social welfare as it lowers the service rate from $\mu$ to $\mu q$. In addition, the socially optimal threshold value is an integer even though customers are allowed to use fractional thresholds. 
Figure \ref{SWNR} (the blue curve) illustrates how the social welfare varies with the threshold value. 

\subsection{Social welfare in the $R$-case}
When customers are allowed to renege after they join, the social welfare calculation is more involved, as not every customer who chooses to join contributes $R_0$ to the social welfare. On this account, to calculate the social welfare in the $R$-case, we need to work out the probability that a joining customer reneges before she successfully completes the service. Denote this probability by $\tilde{p}^{(x)}$ when every customer uses threshold $x$. In order to obtain $\tilde{p}^{(x)}$, we first calculate the distribution 
\[
\tilde{\pi}^{(x)}_{k} = \frac{\mu (1-q) \, \hat{\pi}^{(x)}_{k+1}}{\sum_{k=0}^{\lfloor x \rfloor} \,\mu  (1-q) \, \hat{\pi}^{(x)}_{k+1} } = \frac{\hat{\pi}^{(x)}_{k+1}}{\sum_{k=0}^{\lfloor x \rfloor} \, \hat{\pi}^{(x)}_{k+1} }\qquad k = 0, \cdots, \lfloor x \rfloor \,.
\]
of the number of customers in the system observed by each feedback customer. Each joining customer can only renege when her service fails and there are $\lfloor x \rfloor$ other customers in the system. If this is the case, she reneges with probability $1-p$. So a joining customer reneges at her $k$th feedback with probability
\begin{equation}
\tilde{p}^{(x)}_k : = (1-q)(1-p)\,\tilde{\pi}^{(x)}_{\lfloor x \rfloor} \, \left((1-q) \left(1-(1-p)\,\tilde{\pi}^{(x)}_{\lfloor x \rfloor} \right)\right)^{k-1} \,.
\end{equation}
Hence, the probability that a joining customer reneges before she successfully completes the service is given by 
\begin{align}
\tilde{p}^{(x)} =\sum_{k = 1}^{\infty}\tilde{p}^{(x)}_k
&=\frac{(1-q)\, (1-p)\,\tilde{\pi}^{(x)}_{\lfloor x \rfloor} }{1-(1-q) \, \left(1-(1-p)\,\tilde{\pi}^{(x)}_{\lfloor x \rfloor} \right)} \,.
\end{align}
Then the social welfare in the $R$-case is
\begin{align}
S^{R} (x):&= \lambda \left(\sum_{k = 1}^{\lfloor x \rfloor} \hat{\pi}_{k-1}^{(x)} \, \hat{z}_{k,k}^{(x)} +  (x - \lfloor x \rfloor) \, \hat{\pi}_{\lfloor x \rfloor}^{(x)} \, \hat{z}_{\lfloor x \rfloor+1,\lfloor x \rfloor+1}^{(x)}\right) \, \\
&= \lambda \, R_0 \, \left(\sum_{k = 0}^{\lfloor x \rfloor-1} \hat{\pi}_k^{(x)} +  (x - \lfloor x \rfloor)  \, \hat{\pi}_{\lfloor x \rfloor}^{(x)} \right) \left(1-\tilde{p}^{(x)} \right) -\sum_{k = 0}^{\lceil x \rceil} \, k \, \hat{\pi}_{k}^{(x)} \,.
\end{align}
If we take the derivative of $S^{R}(x)$, we have
\begin{equation} \label{Sod2}
\frac{dS^{R}(x)}{dx} = 
	\Scale[0.9]{
\begin{cases} 
\displaystyle \frac{q \,\rho^{\lfloor x \rfloor} \left(R_0 \lambda (\rho-1)^2 -  \rho \,\left(1 - 2 \rho + \lfloor x \rfloor (1-\rho) + \rho^{\lfloor x \rfloor+2}\right)\right)}{(1 + \rho^{n+1} ((x-\lfloor x \rfloor) \, (1-q \rho)-1))^2}   \qquad & \text{when }x > \lfloor x \rfloor \\
\text{undefined} & \text{when }x = \lfloor x \rfloor \,.
\end{cases}}
\end{equation}
Following a similar argument to that in Proposition \ref{pro1}, there exists a socially optimal threshold $n_R^{\star}$ such that $\dfrac{dS^{R}(x)}{dx}$ is increasing when $ x  \leq n_R^{\star}$; is decreasing when $ x  > n_R^{\star}$. The part in Equation (\ref{Sod2}) that decides the sign of $\displaystyle 	\frac{dS^{R}(x)}{dx}$ is $R_0 \lambda (\rho-1)^2 - \rho \, \left(1 - 2 \rho + \lfloor x \rfloor (1-\rho) + \rho^{\lfloor x \rfloor+2} \right)$, which is the same as in Equation (\ref{Sod1}), so $n_R^{\star} = n_N^{\star}$.

\begin{figure}
	\centering
	\includegraphics[width=10cm]{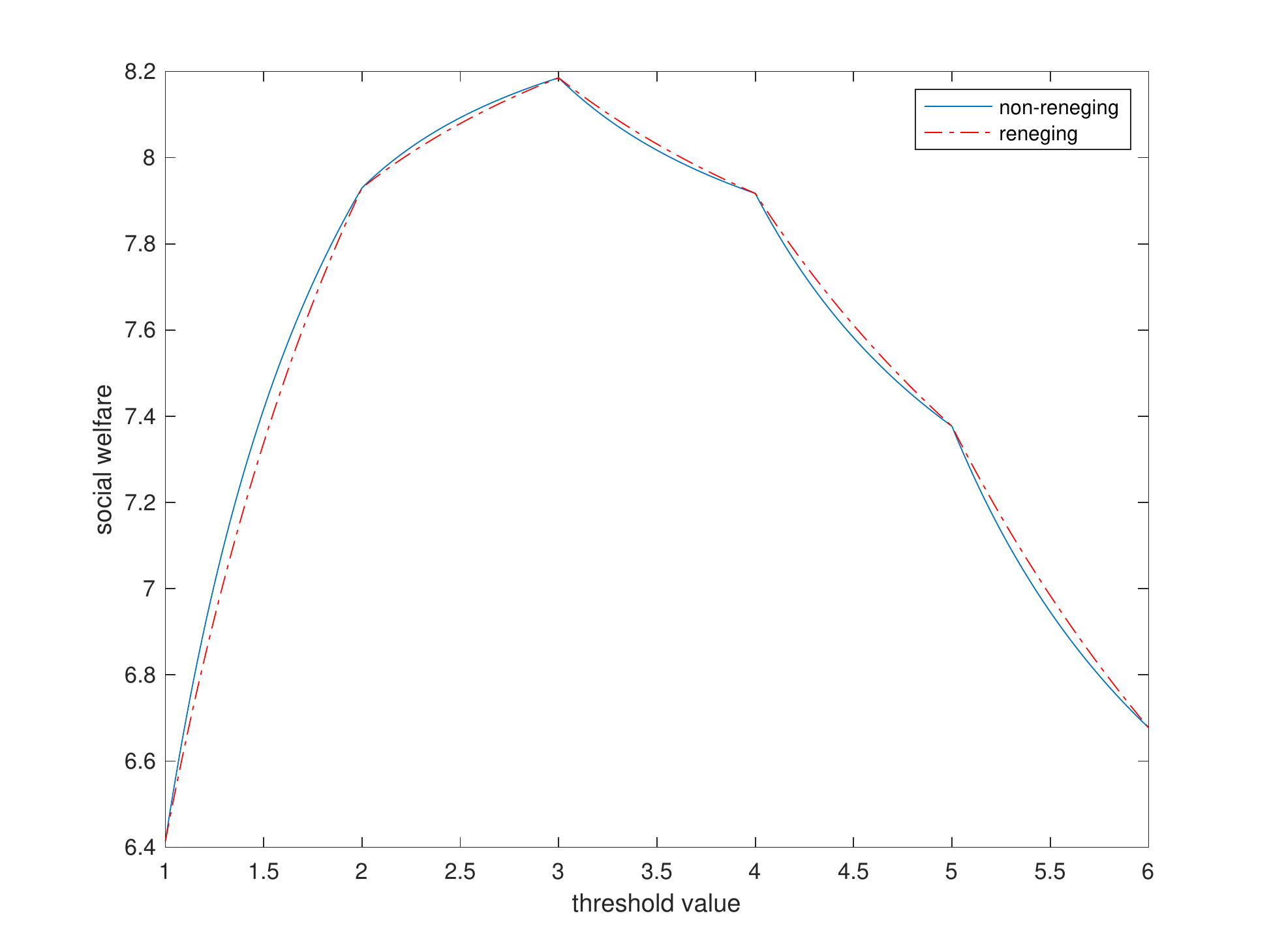}
	\caption{Social welfare when $\lambda = 1, \, \mu = 0.8, \, q = 0.8, \, R_0  = 18$.} \label{SWNR}
\end{figure}

When the threshold is an integer, the joining customers never renege, so there is no difference between the $N$-case and the $R$-case, the socially optimal threshold and the optimal social welfare in both cases are the same. In Figure \ref{SWNR}, an example of social welfare in the $N$-case (blue) and the $R$-case (red) is plotted. The socially optimal threshold $n_N^{\star} = n_R^{\star}=3$. Also, the figure indicates that the social welfare in the non-reneging case is greater than the reneging case when customers use threshold $x < n^*_N$ but is lower when they use $x > n^*_N$. A possible explanation for this is that, when the customers' threshold is less than the socially-desired value, the reason is that fewer customers use the service, and allowing reneging makes this worse. On the other hand, when the customers' threshold is greater than the socially-desired value, the social welfare is less because joining customers inflict negative externalities on others \citep{HO16}, and allowing reneging makes it easier to leave, which improves the situation.


\section*{Acknowledgments}
\noindent P. G. Taylor's research is supported by the Australian Research Council (ARC) Laureate Fellowship FL130100039 and the ARC Centre of Excellence for the Mathematical and Statistical Frontiers (ACEMS). M. Fackrell's research is supported by the ARC Centre of Excellence for the Mathematical and Statistical Frontiers (ACEMS). J. Wang would like to thank the University of Melbourne for supporting her work through the Melbourne Research Scholarship. 

\begin{appendix}
	\normalsize
	\section{The non-reneging case} \label{appendix:1.1}
	\begin{align}
	\begin{split}
	&P^{(x)} = \begin{bmatrix}
	A_0^{(1)} & A_1^{(1)} & 0 & 0 & \cdots & \cdots & \cdots\\
	A_{-1}^{(2)} & A_0^{(2)} & A_1^{(2)}  & 0 & \cdots & \cdots & \cdots\\
	0 & A_{-1}^{(3)} & A_0^{(3)} & A_1^{(3)} & \cdots & \cdots & \cdots \\
	0 & 0 & A_{-1}^{(4)} & A_0^{(4)} & \cdots & \cdots & \cdots \\
	\vdots &\vdots & \vdots & \vdots & \ddots & \ddots &  \vdots \\
	\vdots &\vdots &\vdots & \vdots &  \vdots & A_{-1}^{(\lceil x \rceil + 1)} & A_0^{(\lceil x \rceil + 1)} \\
	\end{bmatrix} \,
	\end{split} \\
	\begin{split}
	&A_{-1}^{(k)} = \begin{bmatrix}
	0 & 0 & \cdots & 0\\
	\frac{\mu q}{\lambda+\mu} & 0 & \cdots & 0 \\
	0 & \frac{\mu q}{\lambda+\mu}  & \cdots & 0 \\
	\vdots & \vdots & \ddots & \vdots  \\
	0 & 0 & \cdots &  \frac{\mu q}{\lambda+\mu}
	\end{bmatrix} \in \mathbb{R}^{k \times (k-1)} \, \qquad k = 2, \cdots, \lceil x \rceil + 1
	\end{split} \\
	\begin{split}
	&A_{0}^{(k)} = 
	\begin{bmatrix}
	0 & 0 & \cdots & \frac{\mu (1-q)}{\lambda+\mu}\\
	\frac{\mu (1-q)}{\lambda+\mu} & 0 & \cdots & 0 \\
	\vdots & \ddots & \ddots & \vdots  \\
	0 & \cdots & \frac{\mu (1-q)}{\lambda+\mu} &  0
	\end{bmatrix}
	\in \mathbb{R}^{k \times k}
	\,  \quad k = 1, \cdots, \lfloor x \rfloor-1
	\end{split}\\
	\begin{split}
	&\Scale[0.95]{A_{0}^{(\lfloor x \rfloor)} = \begin{bmatrix}
	\frac{\lambda(1-(x-\lfloor x \rfloor))}{\lambda+\mu} & 0 & \cdots & 0\\
	0 & \frac{\lambda(1-(x-\lfloor x \rfloor))}{\lambda+\mu} & \cdots & 0 \\
	\vdots & \vdots & \ddots & \vdots  \\
	0 & 0 & \cdots &  \frac{\lambda(1-(x-\lfloor x \rfloor))}{\lambda+\mu}
	\end{bmatrix} +
	\begin{bmatrix}
	0 & 0 & \cdots & \frac{\mu (1-q)}{\lambda+\mu}\\
	\frac{\mu (1-q)}{\lambda+\mu} & 0 & \cdots & 0 \\
	\vdots & \ddots & \ddots & \vdots  \\
	0 & \cdots & \frac{\mu (1-q)}{\lambda+\mu} &  0
	\end{bmatrix}
	\in \mathbb{R}^{\lfloor x \rfloor \times \lfloor x \rfloor}}\, 
	\end{split} \\
	\begin{split}
	&\Scale[0.95]{A_{0}^{(k)} = \begin{bmatrix}
	\frac{\lambda}{\lambda+\mu} & 0 & \cdots & 0\\
	0 & \frac{\lambda}{\lambda+\mu} & \cdots & 0 \\
	\vdots & \vdots & \ddots & \vdots  \\
	0 & 0 & \cdots &  \frac{\lambda}{\lambda+\mu}
	\end{bmatrix} +
	\begin{bmatrix}
	0 & 0 & \cdots & \frac{\mu (1-q)}{\lambda+\mu}\\
	\frac{\mu (1-q)}{\lambda+\mu} & 0 & \cdots & 0 \\
	\vdots & \ddots & \ddots & \vdots  \\
	0 & \cdots & \frac{\mu (1-q)}{\lambda+\mu} &  0
	\end{bmatrix}
	\in \mathbb{R}^{k \times k} \qquad k = \lfloor x \rfloor+1, \, \lceil x \rceil +1}
	\end{split}\\
	\begin{split}
	&A_{1}^{(k)} = 
	\begin{bmatrix}
	\frac{\lambda}{\lambda+\mu} & 0 & \cdots & 0 & 0\\
	0 & \frac{\lambda}{\lambda+\mu} & \cdots & 0 & 0\\
	\vdots & \vdots & \ddots & \vdots & \vdots \\
	0 & 0 & \cdots &  \frac{\lambda}{\lambda+\mu} & 0
	\end{bmatrix}
	\in \mathbb{R}^{k \times (k+1)}
	\,  \quad k = 1, \cdots, \lfloor x \rfloor-1
	\end{split} 
	\end{align}
	\begin{align}
	\begin{split}
	&A_{1}^{(\lfloor x \rfloor)} = \begin{bmatrix}
	\frac{\lambda p}{\lambda+\mu} & 0 & \cdots & 0 & 0\\
	0 & \frac{\lambda p}{\lambda+\mu} & \cdots & 0 & 0\\
	\vdots & \vdots & \ddots & \vdots & \vdots \\
	0 & 0 & \cdots &  \frac{\lambda p}{\lambda+\mu} & 0
	\end{bmatrix}
	\in \mathbb{R}^{\lfloor x \rfloor \times (\lfloor x \rfloor+1)} \qquad A_1^{(\lceil x \rceil)} = \bm{0}_{\lceil x \rceil \times (\lceil x \rceil+1)} \,.
	\end{split}
	\end{align}
	
\section{} \label{appendix:proof of lemma1 and 2}
\subsection{Proof of Lemma \ref{lemmaWI}}
	For integer $d \geq 1$, we first define
	\begin{equation}
	^{d}\nu^{(x)}_{i,j}: = \left( (P^{(x)})^d \bm{e} \right)_{\frac{j(j-1)}{2}+i} \qquad 1 \leq i \leq j \leq \lceil x \rceil \,, 
	\end{equation}
	where $V^d$ is the $d$th power of $V$ and $\bm{v}_k$ is the $k$th element of a  vector $\bm{v}$. We now prove by mathematical induction that $^{d}\nu^{(x)}_{j,j}$ and $^{d}\nu^{(x)}_{i,j}$ are increasing in $j$ for any $d$ and $1 \leq i \leq j \leq \lceil x \rceil$.
	
	When $d =1$,
	\begin{align}
	& ^{1}\nu^{(x)}_{1,j} \,= \,^{1}\nu^{(x)}_{1,j+1} = 1-\frac{\mu q}{\lambda+\mu} \quad 1 \leq j \leq \lceil x \rceil\\
	& ^{1}\nu^{(x)}_{i,j} \,= \,^{1}\nu^{(x)}_{i,j+1} = \,^{1}\nu^{(x)}_{j+1,j+1}= 1 \quad 1 < i \leq j \leq \lceil x \rceil \,.
	\end{align}
	That is, $^{1}\nu^{(x)}_{j+1,j+1} = \, ^{1}\nu^{(x)}_{j,j} > \, ^{1}\nu^{(x)}_{1,1}$ for $1 < j \leq \lceil x \rceil$, and $^{1}\nu^{(x)}_{i,j+1} = \, ^{1}\nu^{(x)}_{i,j}$ for $1 \leq i \leq j \leq \lceil x \rceil$.  
	
	Next suppose that the induction assumption is at the $d$th transition,
	\begin{equation} \label{lemmaWI:eq1}
	^{d}\nu^{(x)}_{j+1,j+1} \, \geq \,  ^{d}\nu^{(x)}_{j,j} \qquad ^{d}\nu^{(x)}_{i,j+1} \, \geq \,  ^{d}\nu^{(x)}_{i,j} \qquad 1 \leq i \leq j \leq \lceil x \rceil \,.
	\end{equation}
	Before proving that \ref{lemmaWI:eq1} holds for $d+1$, we first prove that $^{d}\nu^{(x)}_{i,j} \, \geq \,  ^{d+1}\nu^{(x)}_{i,j}$. Since $^{d}\nu^{(x)}_{i,j}$ represents the sum of probabilities of being in each state in $\mathcal{S}$ at the $d$th transition, that is the probability that the tagged customer is still in the system in the $d$th transition, if the initial state is $(i,j)$. It follows from this physical interpretation that $^{d}\nu_{i,j}^{(x)} = 1$ when $d < i$. Furthermore, since the event that the tagged customer is still in the system after $d+1$ transitions is a subset of the event that it is still in the system after $d$ transitions, it must be the case that $^d\nu^{(x)}_{i,j}$ is decreasing in $d$. 
	
	Then by expanding both $^{d+1}\nu^{(x)}_{i,j+1}$ and $^{d+1}\nu^{(x)}_{i,j}$ as in Equation \eqref{eq:wij} and collecting identical terms together, for $1 \leq i \leq j \leq \lceil x \rceil$, we have
	\begin{align}  
	& \Scale[0.85]{ \left(^{d+1}\nu^{(x)}_{i,j+1}-\, ^{d+1}\nu^{(x)}_{i,j} \right) = \notag }\\
	&\Scale[0.85]{ \frac{\lambda}{\lambda+\mu} \left( \left(^{d}\nu^{(x)}_{i,j+2}-\,^{d}\nu^{(x)}_{i,j+1} \right) \, \mathbbm{1}_{\lbrace j <\lfloor x \rfloor-1 \rbrace} +  p\, \left(^{d}\nu^{(x)}_{i,j+2}-\,^{d}\nu^{(x)}_{i,j+1} \right)  \, \mathbbm{1}_{\lbrace j = \lfloor x \rfloor-1 \rbrace} + (1-p) \left(^{d}\nu^{(x)}_{i,j+1}-\, ^{d}\nu^{(x)}_{i,j} \right) \mathbbm{1}_{\lbrace j = \lfloor x \rfloor \rbrace}  \right)  \notag} \\
	&\Scale[0.85]{ + \frac{\mu }{\lambda+\mu}  \left( (1-q) \, \left(\, ^{d}\nu^{(x)}_{j+1,j+1}-\, ^{d}\nu^{(x)}_{j,j} \right) \, \mathbbm{1}_{\lbrace i =1 \rbrace}  + \left( q \, \left(\,^{d}\nu^{(x)}_{i-1,j} - \, ^{d}\nu^{(x)}_{i-1,j-1} \right)  +  (1-q) \, \left( ^{d}\nu^{(x)}_{i-1,j+1}- \,^{d}\nu^{(x)}_{i-1,j} \right) \right) \mathbbm{1}_{\lbrace  i >1  \rbrace} \right) \label{nu1}} \,,
	\end{align}
	Again, by expanding $^{d+1}\nu^{(x)}_{j+1,j+1}$ using Equation \eqref{eq:wij}, and collapsing the term $\displaystyle \frac{\mu q}{\lambda+\mu} \, ^{d}\nu^{(x)}_{j,j}$, we obtain for $1 \leq i \leq j \leq \lceil x \rceil$,
	\begin{align} 
	& \Scale[0.85]{\left(\, ^{d+1}\nu^{(x)}_{j+1,j+1}- \,^{d+1}\nu^{(x)}_{j,j} \right) \geq \left(^{d+1}\nu^{(x)}_{j+1,j+1}- \, ^{d}\nu^{(x)}_{j,j}\right) = \notag}\\
	& \Scale[0.85]{\frac{\mu (1-q)}{\lambda+\mu}\left(^{d}\nu^{(x)}_{j,j+1}-\,^{d}\nu^{(x)}_{j,j} \right) + \frac{\lambda}{\lambda+\mu} \left( \left(^{d}\nu^{(x)}_{j+1,j+2} - \,^{d}\nu^{(x)}_{j,j} \right) \mathbbm{1}_{\lbrace j < \lfloor x \rfloor-1 \rbrace} \right. \notag} \\
	& \Scale[0.85]{ \left.+ \left(p \, \left(^{d}\nu^{(x)}_{j+1,j+2}- \, ^{d}\nu^{(x)}_{j,j} \right) +(1-p)\, \left(^{d}\nu^{(x)}_{j+1,j+1}- \, ^{d}\nu^{(x)}_{j,j} \right) \right) \mathbbm{1}_{\lbrace j = \lfloor x \rfloor-1 \rbrace} +\left(^{d}\nu^{(x)}_{j+1,j+1}- \,^{d}\nu^{(x)}_{j,j} \right)  \mathbbm{1}_{\lbrace  j = \lfloor x \rfloor  \rbrace}  \label{nu2} \right) \,, }
	\end{align} 
	where the inequality in \eqref{nu2} holds strictly if and only if $d > j-1$.
	It follows from the induction assumption \eqref{lemmaWI:eq1} that $^{d}\nu^{(x)}_{j+1,j+2} \geq  \,^{d}\nu^{(x)}_{j+1,j+1}  \geq \,^{d}\nu^{(x)}_{j,j}$, hence for $1 \leq  i \leq j \leq \lceil x \rceil$,
	\begin{equation}
	( ^{d+1}\nu^{(x)}_{i,j+1} - \, ^{d+1}\nu^{(x)}_{i,j} )\geq 0 \,,
	\end{equation}
	and
	\begin{equation}
	(^{d+1}\nu^{(x)}_{j+1,j+1}- \, ^{d+1}\nu^{(x)}_{j,j}) \geq \frac{\mu(1-q)}{\lambda+\mu} \, (^{d}\nu^{(x)}_{j,j+1}-\, ^{d}\nu^{(x)}_{j,j}) \geq 0 \,.
	\end{equation}

	Finally, from equation (\ref{eq: poisson1})
	\begin{equation}
	\bm{w}^{(x)} = \frac{1}{\lambda+\mu} \, (I-P^{(x)})^{-1} \bm{e} =\frac{1}{\lambda+\mu} \, \sum_{d=0}^{\infty} (P^{(x)})^d \bm{e} \,.
	\end{equation}
	It follows that
	\begin{align}
	& w_{i,j} = \frac{1}{\lambda+\mu} \left(1+ \sum_{d = 1}^{\infty} \, ^{d}\nu^{(x)}_{i,j} \right) \,, \label{lemma1_eq3}\\
	& w_{j,j} = \frac{1}{\lambda+\mu} \left(1+ \sum_{d = 1}^{\infty} \, ^{d}\nu^{(x)}_{j,j} \right) 
	\end{align}
	are increasing in $j$. 
	\hfill $\square$
	
\subsection{Proof of Lemma \ref{lemmaWx}}
\begin{itemize}
	\item When $n <  x_1 = n+ p_1 < n+ p_2 = x_2 \leq n+1$, from equation (\ref{eq: poisson1}),
	\begin{align}
	& \left(\mathbf{I}-P^{(x_1)} \right)
	\bm{w}^{(x_1)} \, 
	=  \, \frac{1}{\lambda+\mu}\bm{e}\\
	& \left(\mathbf{I}-P^{(x_2)} \right)
	\bm{w}^{(x_2)} \, 
	=  \, \frac{1}{\lambda+\mu}\bm{e} \,,
	\end{align} 
	where the probability transition matrices $P^{(x_1)}$ and $P^{(x_2)}$ have the same dimension. Noting that $P^{(x_2)}$ and $P^{(x_1)}$ only differ in $n$ rows, we have
	\begin{equation} \label{lemmaWx:eq1}
	\left(I-P^{(x_2)} \right) \, (\bm{w}^{(x_2)}-\bm{w}^{(x_1)}) = (P^{(x_2)}-P^{(x_1)}) \bm{w}^{(x_1)} =
	\frac{\lambda(p_2-p_1)}{\lambda+\mu} \,  
	\begin{bmatrix}
	\mathbf{0}_{\frac{n(n-1)}{2} \times 1} \\
	w^{(x_1)}_{1,n+1} - w^{(x_1)}_{1,n} \\
	w^{(x_1)}_{2,n+1} - w^{(x_1)}_{2,n} \\
	\vdots \\
	w^{(x_1)}_{n,n+1} - w^{(x_1)}_{n,n} \\
	\mathbf{0}_{(2n+3) \times 1}
	\end{bmatrix} \,.
	\end{equation}
	We know from equation (\ref{lemma1_eq3}) that $w^{(x_1)}_{i,n+1} - w^{(x_1)}_{i,n} > 0$.
	Also, since the $((i_1,j_1),(i_2,j_2))$th entry of  $(I-P^{(x_2)})^{-1}$ is the tagged customer's expected number of visits to state $(i_2,j_2)$ starting from state $(i_1,j_1)$ in $\mathcal{S}$ before she leaves the system, 
	\[
	(I-P^{(x_2)})^{-1} > 0 \,.
	\]
	Hence
	\begin{equation}  \label{lemmaWx:eq3}
	\bm{w}^{(x_2)}-\bm{w}^{(x_1)} = (I-P^{(x_2)})^{-1} \, (P^{(x_2)}-P^{(x_1)}) \, \bm{w}^{(x_1)} > 0 \,. 
	\end{equation}	
	
	\item When $x_1 = n$ and $x_2 =n + p_2 $ or $x_2 = n+1$, $P^{(x_1)}$ and $P^{(x_2)}$ have different sizes. However, we can write $P^{(x_2)}$ as
	\begin{align}
	P^{(x_2)} \quad &= \quad \left[\begin{array}{c c c c c | c}
	A_0^{(1)} & A_1^{(1)} & \cdots & \cdots & \cdots  & 0\\
	A_{-1}^{(2)} & A_0^{(2)} & A_1^{(2)}  & \cdots & \cdots & 0\\
	\vdots & A_{-1}^{(3)} & A_0^{(3)} & A_1^{(3)}   & \cdots & 0 \\
	\vdots &\vdots & \vdots & \ddots  & \ddots  & \vdots \\
	\vdots & \vdots & \vdots & A_{-1}^{(n+1)} & A_0^{(n+1)} &  0 \\
	\hline
	0 & \cdots & \cdots & \cdots & {A}_{-1}^{(n+2)} & {A}_0^{(n+2)} 
	\end{array} \right] \\
	&= \quad\left[\begin{array}{c  | c}
	\bar{P}^{(x_2)} & \bm{0}_{\frac{(n+1)(n+2)}{2} \times (n+2)}\\
	\hline
	\bm{0}_{(n+2) \times \frac{n(n+1)}{2}} \quad {A}_{-1}^{(n+2)}& {A}_0^{(n+2)} 
	\end{array} \right] \,,
	\end{align}
	and
	\begin{equation}
	(I-P^{(x_2)})^{-1} = \left[\begin{array}{c | c}
	( I-\bar{P}^{(x_2)})^{-1}    & \mathbf{0} \\
	\hline
	N \, ( I-\tilde{P}^{(x_2)})^{-1}
	& ( I- {A}_0^{(n+2)}  )^{-1}  
	\end{array} \right] \,,
	\end{equation}
	where $N = \begin{bmatrix}
	\bm{0}_{(n+2) \times \frac{n(n+1)}{2}} & ( I-  {A}_0^{(n+2)}  )^{-1} \, {A}_{-1}^{(n+2)}
	\end{bmatrix}$. When $x_1 = n$, the position where the tagged customer can join is at most $n+1$, so we compare $\bm{w}^{(x_1)}_k$ with $\bm{w}^{(x_2)}_k$ for $k = 1, \cdots, \displaystyle \frac{(n+1)(n+2)}{2}$. If we define $\bm{\bar{w}}^{(x_2)}: = \left[{I}_{\frac{(n+1)(n+2)}{2}} \quad \bm{0}\right] \bm{w}^{(x_2)} $, then
	\begin{equation}
	\bm{\bar{w}}^{(x_2)} = \frac{1}{\lambda+\mu}( I-\bar{P}^{(x_2)})^{-1} \, \bm{e} \,,
	\end{equation}
	thus
	\begin{equation} \label{lemmaWx:eq2}
	\bm{\bar{w}}^{(x_2)}-\bm{w}^{(x_1)} = \frac{\lambda \, p_2}{\lambda+\mu} \,  (I-\bar{P}^{(x_2)})^{-1}\, 
	\begin{bmatrix}
	\mathbf{0}_{\frac{n(n-1)}{2} \times 1} \\
	w^{(x_1)}_{1,n+1} - w^{(x_1)}_{1,n} \\
	w^{(x_1)}_{2,n+1} - w^{(x_1)}_{2,n} \\
	\vdots \\
	w^{(x_1)}_{n,n+1} - w^{(x_1)}_{n,n} \\
	\mathbf{0}_{(n+1) \times 1}
	\end{bmatrix} > 0 \,.
	\end{equation}
	In \eqref{lemmaWx:eq2}, with an abuse of notation, we include the case $p_2 = 1$.
	\item When $x_1 = n_1 +p_1, x_2= n_2 +p_2$, and $n_1 <n_2$, by comparing the expected sojourn time for all the consecutive integers between $x_1$ and $x_2$, it follows from the aforementioned reasoning that 
	\begin{align}
	& w_{i,j}^{(x_1)} < w_{i,j}^{(\lceil x_1 \rceil)}   & 1 & \leq i \leq j \leq \lceil x_1 \rceil +1 \\
	& \cdots &  \notag\\
	& w_{i,j}^{(\lceil x_2 \rceil-2)} < w_{i,j}^{(\lceil x_2 \rceil-1)}   & 1 &\leq i \leq j \leq \lceil x_2 \rceil -1 \\
	& w_{i,j}^{(\lceil x_2 \rceil-1)} < w_{i,j}^{(x_2)}   & 1 &\leq i \leq j \leq \lceil x_2 \rceil \,.
	\end{align}
Hence $ w_{i,j}^{(x_1)} <  w_{i,j}^{(x_2)}, \, 1 \leq i \leq j \leq \lceil x_1\rceil +1$.
\end{itemize}
$\hfill \square$

\section{The reneging case}
\subsection{} \label{appendix: 2.1}
\begin{equation}
\hat{P}^{(\lfloor x \rfloor+1,x)} = \begin{bmatrix}
A_0^{(1)} & A_1^{(1)} & \cdots & \cdots & \cdots & \cdots & 0\\
A_{-1}^{(2)} & A_0^{(2)} & A_1^{(2)}  & \cdots & \cdots & \cdots & 0\\
\vdots & A_{-1}^{(3)} & A_0^{(3)} & A_1^{(3)}   & \cdots & \cdots & 0 \\
\vdots & \vdots & A_{-1}^{(4)} & A_0^{(4)} & A_1^{(4)} & \cdots & 0 \\
\vdots &\vdots & \vdots & \vdots & \vdots & \ddots & \vdots \\
0 &0 &0 & 0 &  0 & \hat{A}_{-1}^{(\lfloor x \rfloor+1)} & \tilde{A}_0^{(\lfloor x \rfloor+1)} 
\end{bmatrix} \,,
\end{equation}
\begin{align}
&\Scale[0.95]{\hat{A}_{-1}^{(\lfloor x \rfloor+1)} = \begin{bmatrix}
		0 & 0 & \cdots & 0\\
		\frac{\mu q + \mu (1-q)(1-(x-\lfloor x \rfloor))}{\lambda+\mu} & 0 & \cdots & 0 \\
		0 & \frac{\mu q + \mu (1-q)(1-(x-\lfloor x \rfloor))}{\lambda+\mu}  & \cdots & 0 \\
		\vdots & \vdots & \ddots & \vdots  \\
		0 & 0 & \cdots &  \frac{\mu q + \mu (1-q)(1-(x-\lfloor x \rfloor))}{\lambda+\mu}
		\end{bmatrix} \in \mathbb{R}^{(\lfloor x \rfloor + 1) \times  \lfloor x \rfloor}} \\
&\Scale[0.95]{ \tilde{A}_0^{(\lfloor x \rfloor+1)} =  \begin{bmatrix}
		\frac{\lambda}{\lambda+\mu} & 0 & \cdots & 0\\
		0 & \frac{\lambda}{\lambda+\mu} & \cdots & 0 \\
		\vdots & \vdots & \ddots & \vdots  \\
		0 & 0 & \cdots &  \frac{\lambda}{\lambda+\mu}
		\end{bmatrix} +
		\begin{bmatrix}
		0 & 0 & \cdots & \frac{\mu (1-q)}{\lambda+\mu}
		\\
		\frac{\mu (1-q)\, (x - \lfloor x \rfloor)}{\lambda+\mu} & 0 & \cdots & 0 \\
		\vdots & \ddots & \ddots & \vdots  \\
		0 & \cdots & \frac{\mu (1-q) \, (x - \lfloor x \rfloor)}{\lambda+\mu} &  0
		\end{bmatrix}
		\in \mathbb{R}^{(\lfloor x \rfloor+1)  \times (\lfloor x \rfloor+1 )} }\\
& \bm{g} = \frac{\mu q  \, R_0}{\lambda+\mu}  \, \bm{{e}}_k \, - \, \frac{1}{\lambda+\mu} \, \bm{{e}} \qquad k = \frac{j(j-1)}{2}+1, \, j = 1, \ldots, \lfloor x \rfloor+1 \,.
\end{align}
\subsection{} \label{appendix: 3.1}
\begin{align}
&\hat{P}^{(x,x)} = \begin{bmatrix}
A_0^{(1)} & A_1^{(1)} & \cdots & \cdots & \cdots & \cdots & 0\\
A_{-1}^{(2)} & A_0^{(2)} & A_1^{(2)}  & \cdots & \cdots & \cdots & 0\\
\vdots & A_{-1}^{(3)} & A_0^{(3)} & A_1^{(3)}   & \cdots & \cdots & 0 \\
\vdots & \vdots & A_{-1}^{(4)} & A_0^{(4)} & A_1^{(4)} & \cdots & 0 \\
\vdots &\vdots & \vdots & \vdots & \vdots & \ddots & \vdots \\
0 &0 &0 & 0 &  0 & \hat{A}_{-1}^{(\lfloor x \rfloor+1)} & \hat{A}_0^{(\lfloor x \rfloor+1)} 
\end{bmatrix} \,, \\
&\Scale[0.95]{\hat{A}_{0}^{(\lfloor x \rfloor+1)} = \begin{bmatrix}
		\frac{\lambda}{\lambda+\mu} & 0 & \cdots & 0\\
		0 & \frac{\lambda}{\lambda+\mu} & \cdots & 0 \\
		\vdots & \vdots & \ddots & \vdots  \\
		0 & 0 & \cdots &  \frac{\lambda}{\lambda+\mu}
		\end{bmatrix} +
		\begin{bmatrix}
		0 & 0 & \cdots & \frac{\mu (1-q) (x-\lfloor x \rfloor)}{\lambda+\mu}
		\\
		\frac{\mu (1-q) (x-\lfloor x \rfloor)}{\lambda+\mu} & 0 & \cdots & 0 \\
		\vdots & \ddots & \ddots & \vdots  \\
		0 & \cdots & \frac{\mu (1-q) (x-\lfloor x \rfloor)}{\lambda+\mu} &  0
		\end{bmatrix}
		\in \mathbb{R}^{(\lfloor x \rfloor+1) \times (\lfloor x \rfloor+1)} }
\end{align}

\section{Proof of Lemma \ref{lemma:CompareRNR1}} \label{appendix:lemma3}
Proof. It follows from the explanation at the beginning of this section that the tagged customer will never renege after joining if she uses $x_{tag} \geq \lceil x \rceil$ when others uses $x$. Thus the comparison of $z^{(x)}_{i,j}$ and $\hat{z}_{i,j}^{(\lfloor x \rfloor +1,x)}$ is actually the comparison of the respective expected sojourn time $w^{(x)}_{i,j}$ and $\hat{w}_{i,j}^{(\lfloor x \rfloor +1,x)}$. As in Lemma \ref{lemmaWI}, the tagged customer's expected sojourn time is her expected total number of visits of different states until she leaves the system, times $\frac{1}{\lambda+\mu}$. Similar to the proof of Lemma \ref{lemmaWI}, we first define
\begin{equation}
^{d}\hat{\nu}^{(x)}_{i,j}: = \left( (\hat{P}^{(\lfloor x \rfloor+1, x)})^d \bm{e} \right)_{\frac{j(j-1)}{2}+i} \qquad 1 \leq i \leq j \leq \lfloor x \rfloor +1 \,,
\end{equation}
which is the probability that the tagged customer is still in the system at the $d$th transition, if her initial state is $(i,j)$, she never reneges, and others join and renege with threshold $x$. We prove Equation \eqref{eq:compare1} by mathematical induction. 

\begin{figure} 
	\subcaptionbox{start with state $(2,3)$}%
	{\includegraphics[width=0.48\linewidth]{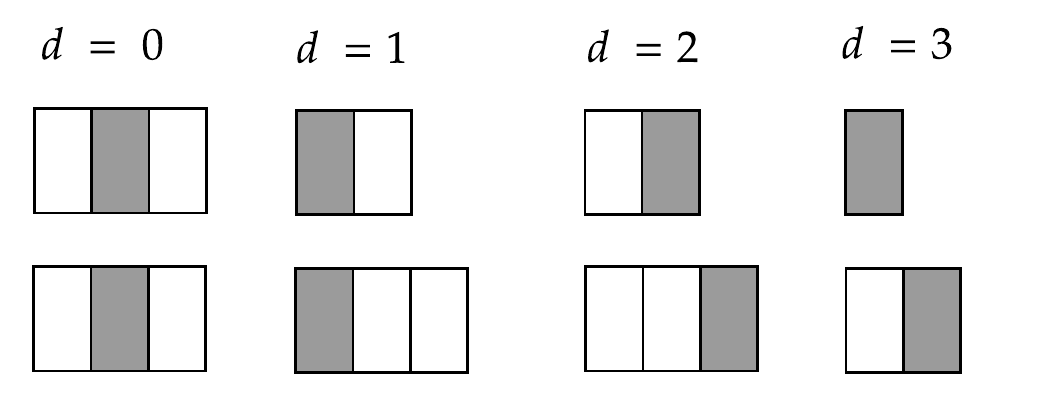}}  
	\hspace{\fill}
	\subcaptionbox{start with state $(3,3)$}%
	{\includegraphics[width=0.46\linewidth]{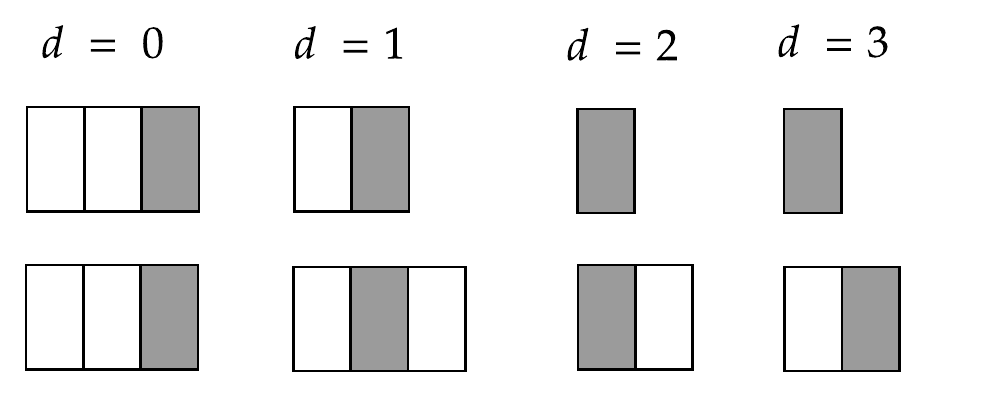}} 
	\caption{An example of the comparison between $^{d} \nu^{(x)}_{i,\lfloor x \rfloor +1}$ and $^{d}\hat{\nu}^{(\lfloor x \rfloor+1,x)}_{i,\lfloor x \rfloor +1}$ when $2 < x < 3$.}   \label{fig:RNR}
\end{figure}

First, it follows from the physical interpretation that when $d \leq \lfloor x \rfloor$, $^d \nu^{(x)}_{i,j} = \, ^d\hat{\nu}^{(\lfloor x \rfloor+1,x)}_{i,j}$. This is because others' reneging will not affect the tagged customer until she rejoins the queue and is in service for the second time, and the queue size has reached $\lfloor x \rfloor + 1$ before she is in service for the first time. This is possible only when $d > \lfloor x \rfloor$. Indeed,
\begin{equation} 
\label{eq:RNR1}
^{\lfloor x \rfloor+1} \nu^{(x)}_{i,j} > \, ^{\lfloor x \rfloor+1}\hat{\nu}^{(\lfloor x \rfloor+1,x)}_{i,j} \qquad1< i \leq j = \lfloor x \rfloor + 1 \,.
\end{equation}
Figure \ref{fig:RNR} shows an example of Equation \eqref{eq:RNR1} with $(i,j) = (2,3)$ in (a), $(i,j) = (3,3)$ in (b), and $2 < x < 3$. The gray block represents the tagged customer and the white ones represent other customers in the system. It shows the sample paths where the reneging of others affects the tagged customer at the earliest possible transition. The first row is when others use threshold $2$, and the second row is when others use threshold $3$. In both Figure \ref{fig:RNR}(a) and (b), the reneging of others affects the queue size at the first transition, but will not affects the tagged customer until $d >3$. 

Next, we assume that $ ^d \nu^{(x)}_{i,j} \geq \,^d\hat{\nu}^{(\lfloor x \rfloor+1,x)}_{i,j}$. Then we write the difference between $P^{(x)}$ and $\hat{P}^{(\lfloor x \rfloor +1,x)}$ in the form
\begin{equation} \label{NRC}
P^{(x)}= \hat{P}^{(\lfloor x \rfloor +1,x)}    + E \Delta \,,
\end{equation}
where
\begin{equation}
	E = \begin{bmatrix}
		\mathbf{0}_{\left(\frac{\lfloor x \rfloor (\lfloor x \rfloor +1)}{2}+1\right) \times \lfloor x \rfloor } \\
		{I}_n
		\end{bmatrix} \quad
		\Delta = \begin{bmatrix}
		\mathbf{0}_{\lfloor x \rfloor \times  \frac{\lfloor x \rfloor (\lfloor x \rfloor -1)}{2}} & -\frac{\mu(1-q)(1-p)}{\lambda+\mu} {I}_{\lfloor x \rfloor}  & \frac{\mu(1-q)(1-p)}{\lambda+\mu} {I}_{\lfloor x \rfloor } & \mathbf{0}_{\lfloor x \rfloor  \times 1} 
		\end{bmatrix} \,.
\end{equation}
Hence,
\begin{equation}
(P^{(x)})^{d+1} \bm{e} = P^{(\lfloor x \rfloor+1,x)}\,\left( (P^{(x)})^d \bm{e} \right)  + E \Delta\,\left( (P^{(x)})^d \bm{e} \right)  \geq   P^{(\lfloor x \rfloor+1,x)} \,\left( (P^{(x)})^d \bm{e} \right) \geq (P^{(\lfloor x \rfloor+1,x)})^{d+1} \bm{e}\,.
\end{equation}
The first inequality follows from the conclusion that
$^{d}\nu^{(x)}_{i,j}$ is increasing in $j$ for any $d$ and $1 \leq i \leq j \leq \lceil x \rceil$ in Lemma \ref{lemmaWI}. Since $^{d}\nu^{(x)}_{i,\lfloor x \rfloor + 1} \geq ^{d}\nu^{(x)}_{i,\lfloor x \rfloor}$ for $i = 1, \dots, \lfloor x \rfloor$, 
\begin{equation}
\Delta\,\left( (P^{(x)})^d \bm{e} \right) = \frac{\mu(1-q)(1-p)}{\lambda+\mu}
\begin{bmatrix}
&^{d}\nu^{(x)}_{1,\lfloor x \rfloor + 1}  - ^{d}\nu^{(x)}_{1,\lfloor x \rfloor} \\
& ^{d}\nu^{(x)}_{2,\lfloor x \rfloor + 1}  - ^{d}\nu^{(x)}_{2,\lfloor x \rfloor} \\
& \vdots \\
& ^{d}\nu^{(x)}_{\lfloor x \rfloor,\lfloor x \rfloor + 1}  - ^{d}\nu^{(x)}_{\lfloor x \rfloor,\lfloor x \rfloor} 
\end{bmatrix} \geq 0 \,.
\end{equation}
The second inequality is from the
induction assumption. Hence,
\begin{equation}
^{d+1} \nu^{(x)}_{i,j} \geq \, ^{d+1} \hat{\nu}^{(\lfloor x \rfloor+1,x)}_{i,j} \,.
\end{equation} 
This concludes the proof for Inequality \eqref{eq:compare1}.

When $x = \lfloor x \rfloor$, $\hat{P}^{(\lfloor x \rfloor, \lfloor x \rfloor)} = P^{(\lfloor x \rfloor)}$, and so $\bm{\hat{z}}^{(\lfloor x \rfloor,\lfloor x\rfloor)} = \bm{z}^{(\lfloor x \rfloor)}$. Also, $A_1^{(\lfloor x \rfloor)} = 0$, so
\begin{align} \label{eq: Px}
& \left[ I_{\frac{\lfloor x \rfloor (\lfloor x \rfloor+1)}{2}}  \quad 0_{\frac{\lfloor x \rfloor (\lfloor x \rfloor+1)}{2} \times (\lfloor x \rfloor +1)} \right] P^{(\lfloor x \rfloor )} \\
=&
{\left[ I_{\frac{\lfloor x \rfloor (\lfloor x \rfloor+1)}{2}}  \quad 0_{\frac{\lfloor x \rfloor (\lfloor x \rfloor+1)}{2} \times (\lfloor x \rfloor +1)} \right] P^{(\lfloor x \rfloor )} 
		\begin{bmatrix}
		I_{\frac{\lfloor x \rfloor (\lfloor x \rfloor+1)}{2}}  \notag \\
		0_{\frac{\lfloor x \rfloor (\lfloor x \rfloor+1)}{2} \times (\lfloor x \rfloor +1)}
		\end{bmatrix}
		\left[ I_{\frac{\lfloor x \rfloor (\lfloor x \rfloor+1)}{2}}  \quad 0_{\frac{\lfloor x \rfloor (\lfloor x \rfloor+1)}{2} \times (\lfloor x \rfloor +1)} \right] \,. }
\end{align}
Multiplying equation (\ref{eq: poisson1}) with $\left[ I_{\frac{\lfloor x \rfloor (\lfloor x \rfloor+1)}{2}}  \quad 0_{\frac{\lfloor x \rfloor (\lfloor x \rfloor+1)}{2} \times (\lfloor x \rfloor +1)} \right]$ on both sides and applying equation (\ref{eq: Px}), we have
\begin{equation}
(I - \bar{P}^{(\lfloor x \rfloor )}) \, \left[ I_{\frac{\lfloor x \rfloor (\lfloor x \rfloor+1)}{2}}  \quad 0_{\frac{\lfloor x \rfloor (\lfloor x \rfloor+1)}{2} \times (\lfloor x \rfloor +1)} \right] \bm{w}^{(\lfloor x \rfloor )} = \frac{1}{\lambda+\mu} \, \bm{e}
\end{equation}
where $\bar{P}^{(\lfloor x \rfloor )}:= \left[ I_{\frac{\lfloor x \rfloor (\lfloor x \rfloor+1)}{2}}  \quad 0_{\frac{\lfloor x \rfloor (\lfloor x \rfloor+1)}{2} \times (\lfloor x \rfloor +1)} \right] P^{(\lfloor x \rfloor )} 
\begin{bmatrix}
I_{\frac{\lfloor x \rfloor (\lfloor x \rfloor+1)}{2}}  \notag \\
0_{\frac{\lfloor x \rfloor (\lfloor x \rfloor+1)}{2} \times (\lfloor x \rfloor +1)}
\end{bmatrix}$. Hence,
\begin{align}
(\bm{w}^{(\lfloor x \rfloor )})_{1:\frac{\lfloor x \rfloor (\lfloor x \rfloor +1)}{2}} = \left( I - \bar{P}^{(\lfloor x \rfloor )} \right)^{-1}\frac{1}{\lambda+\mu} \, \bm{e} \,.
\end{align}
Similarly,
\begin{equation}
(\bm{\hat{w}}^{(\lfloor x \rfloor +1, \lfloor x \rfloor)})_{1:\frac{\lfloor x \rfloor (\lfloor x \rfloor +1)}{2}} = \left( I - \bar{\hat{P}}^{(\lfloor x \rfloor +1, \lfloor x \rfloor )} \right)^{-1}\frac{1}{\lambda+\mu} \, \bm{e} \,,
\end{equation}
where $\bar{\hat{P}}^{(\lfloor x \rfloor )}:= \left[ I_{\frac{\lfloor x \rfloor (\lfloor x \rfloor+1)}{2}}  \quad 0_{\frac{\lfloor x \rfloor (\lfloor x \rfloor+1)}{2} \times (\lfloor x \rfloor +1)} \right] \hat{P}^{(\lfloor x \rfloor +1, \lfloor x \rfloor )}
\begin{bmatrix}
I_{\frac{\lfloor x \rfloor (\lfloor x \rfloor+1)}{2}}  \notag \\
0_{\frac{\lfloor x \rfloor (\lfloor x \rfloor+1)}{2} \times (\lfloor x \rfloor +1)}
\end{bmatrix}$.
Since $P^{(\lfloor x \rfloor+1,  x  )}$ and $P^{(x)}$ have the same first $\displaystyle \frac{\lfloor x \rfloor (\lfloor x \rfloor +1)}{2}$ lines and rows, it follows that $\bm{w}_{1:\frac{\lfloor x \rfloor (\lfloor x \rfloor +1)}{2}}^{(\lfloor x \rfloor )} = \bm{\hat{w}}_{1:\frac{\lfloor x \rfloor (\lfloor x \rfloor +1)}{2}}^{(\lfloor x \rfloor +1, \lfloor x \rfloor )}$. 
\hfill $\square$ 

\section{}
\subsection{The derivative of function $f(p)$} \label{appendix:fp} 
\begin{align}
	\frac{d f(p)}{d p} \notag & =  -\dfrac{\splitfrac{\lambda  \left(2 \lambda ^2 p^2 q (\lambda +\mu  q)+2 \lambda  p (\lambda +\mu  q) (\lambda (1-q)+\mu  q ((q-3) q+3)) \right.}{ \left.+\mu  (1-q) \left(2 \lambda ^2+\mu ^2 (2-q) q^2 ((q-3) q+3)+\lambda  \mu  q ((q-5) q+8)\right)\right)}}{(\mu +\lambda  p)^2 (\lambda  p (\lambda +\mu  q)+\mu  (\lambda +\mu  q ((q-3) q+3)))^2}\\
	&< 0
\end{align}
The calculation of $f'(p)$ is implemented in 
Wolfram Mathematica. 
\subsection{Social welfare expressions} \label{appendix:sw}
The social welfare in the $N-$case
\begin{align}
	& S^{N} (x) = \notag \\
	&  \Scale[1]{\dfrac{\splitfrac{\lambda R_0 (\rho - 1) \, ((1+(x-\lfloor x \rfloor) \, (\rho-1)) \, \rho^{\lfloor x \rfloor} \, - 1) }{+ \rho \, ((1 - \lfloor x \rfloor (1 + (x-\lfloor x \rfloor) (\rho - 1)) ( \rho-1) - (x-\lfloor x \rfloor) (\rho-1)^2) \rho^{\lfloor x \rfloor} - 1)}}{1 + \rho ((1 + (x-\lfloor x \rfloor) (\rho - 1)) (\rho - 1) \rho^{\lfloor x \rfloor} - 1)} }\,.
\end{align}
The social welfare in the $R-$case
\begin{align}
	&S^{R} (x) = \notag \\
	& \Scale[0.9]{\dfrac{\splitfrac{\rho  \left( R_0 \, \mu  q (1-\rho) \left(\rho ^n ((x-\lfloor x \rfloor) (q \rho -1)+1)- (x-\lfloor x \rfloor) (1-q)+1\right) \right.}{+ \left.\left(\rho ^n ((n (\rho -1) ((x-\lfloor x \rfloor) ( \rho -1)+1)+(x-\lfloor x \rfloor) (q (\rho -2) \rho +1)-1))- (x-\lfloor x \rfloor) (1-q)+1\right)\right)}}{(1-\rho) \left(\rho ^{n+1} ((x-\lfloor x \rfloor)  (q \rho -1)+1)+ (x-\lfloor x \rfloor) (1-q)-1\right)}} \,.
\end{align}
\end{appendix}

\end{document}